\documentclass[11pt]{article}
\usepackage{amsmath,amssymb,color,url,mathtools,amsthm}
\usepackage[mathscr]{euscript}
\usepackage{float}
\usepackage{graphicx}
\usepackage[left=1.2in,right=1.2in,top=1in,bottom=1in]{geometry}
\usepackage{tikz}
\usepackage{tikz-cd}
\usepackage{epstopdf}

\usepackage[numbers]{natbib}

\newcounter{commentcounter}

\newtheorem{thm}{Theorem}[section]
\newtheorem*{thm*}{Theorem}
\usepackage{epsfig}
\newtheorem{lem}[thm]{Lemma}
\newtheorem{prop}[thm]{Proposition}
\newtheorem{cor}[thm]{Corollary}
\theoremstyle{definition}
\newtheorem{definition}[thm]{Definition}
\newtheorem{example}[thm]{Example}

\theoremstyle{remark}
\newtheorem*{remark}{Remark}
\theoremstyle{comments}

\theoremstyle{notation}
\newtheorem*{notation}{Notation}
\newtheorem*{mainquestion}{Main Question}
\newtheorem{conjecture}{Conjecture}
\newtheorem{observation}{Observation}

\newtheoremstyle{named}{}{}{\itshape}{}{\bfseries}{.}{.5em}{\thmnote{#3}}
\theoremstyle{named}
\newtheorem*{thmn}{Theorem}

\author{Michael Harrison}

\title{Introducing totally nonparallel immersions}

\newcommand{\Proj}{\mathbb P}
\newcommand{\R}{\mathbb R}
\newcommand{\C}{\mathbb C}
\newcommand{\N}{\mathbb N}

\newcommand{\Z}{\mathbb Z}

\newcommand{\TN}{\textnormal{TN}}
\newcommand{\TS}{\textnormal{TS}}

\newcommand{\Ker}{\textnormal{Ker}}
\newcommand{\rank}{\textnormal{rank}}
\newcommand{\Coker}{\textnormal{Coker}}

\newcommand{\Hom}{\textnormal{Hom}}

\newcommand{\RP}{\mathbb R \mathbb P}

\newcommand{\codim}{\textnormal{codim}}
\newcommand{\dist}{\operatorname{dist}}
\newcommand{\himm}{\mathscr{H}_{\mbox{imm}}}
\newcommand{\hman}{\mathscr{H}_{\mbox{man}}}

\newcommand{\hnot}{\mathscr{H}_{\mbox{not}}}

\begin{document}

\maketitle

\begin{abstract}
An immersion of a smooth $n$-dimensional manifold $M \to \R^q$ is called \emph{totally nonparallel} if, for every distinct $x, y \in M$, the tangent spaces at $f(x)$ and $f(y)$ contain no parallel lines; equivalently, they span a $2n$-dimensional space.  Given a manifold $M$, we seek the minimum dimension $\TN(M)$ such that there exists a totally nonparallel immersion $M \to \R^{\TN(M)}$.  In analogy with the \emph{totally skew embeddings} studied by Ghomi and Tabachnikov, we find that totally nonparallel immersions are related to the generalized vector field problem, the immersion and embedding problems for real projective spaces, and nonsingular symmetric bilinear maps.

Our study of totally nonparallel immersions follows a recent trend of studying conditions which manifest on the \emph{configuration space} $F_k(M)$ of $k$-tuples of distinct points of $M$; for example, $k$-regular embeddings, $k$-skew embeddings, $k$-neighborly embeddings, and several others.  Typically, a map satisfying one of these \emph{configuration space conditions} induces some $S_k$-equivariant map on the configuration space $F_k(M)$ (or on a bundle thereof) and obstructions can be computed in the form of Stiefel-Whitney classes.  However, the existence problem for such conditions is relatively unstudied.

Our main result is a Whitney-type theorem: every smooth $n$-manifold $M$ admits a totally nonparallel immersion into $\R^{4n-1}$, one dimension less than given by genericity.  We begin by studying the local problem, which requires a thorough understanding of the space of nonsingular symmetric bilinear maps, after which the main theorem is established using the removal-of-singularities h-principle technique due to Gromov and Eliashberg.  When combined with a recent non-immersion theorem of Davis, we obtain the exact value $\TN(\RP^n) = 4n-1$ when $n$ is a power of $2$.  This is the first optimal-dimension result for any closed manifold $M$ besides $S^1$, for any of the recently-studied configuration space conditions.
\end{abstract}

\section{Introduction and Main Results}
\label{sec:intro}

\renewcommand{\thefootnote}{\fnsymbol{footnote}} 
\footnotetext{\emph{2010 Mathematics Subject Classification} 57R45 (Primary), 58A20, 58K20 (Secondary)}

In a 1968 paper \cite{Segre}, B.\ Segre constructed a \emph{skew loop} -- a loop in $\R^3$ without parallel tangent lines -- thereby disproving a short-lived nonexistence conjecture of Steinhaus.  Skew loops were used by Ghomi in his study of the ``shadow problem'' \cite{GhomiShadow} and further popularized by Ghomi and Solomon \cite{GhomiSolomon}, who showed that skew loops cannot live on quadric surfaces.  An explicit example of a skew loop may be found in \cite{GhomiShadow} and is depicted in Figure \ref{fig:skewloop} below.  The study of skew loops continued as Solomon constructed skew loops in flat tori \cite{Solomon}, Wu constructed skew loops in every knot class \cite{Wu}, and Ghomi studied skew loops in higher-dimensional spaces and their relationship to ellipsoids \cite{GhomiEllipsoids}.

\begin{figure}[h!t]
    \centering \ \ \ \ \ \ \
    \begin{minipage}{0.36\textwidth}
        \centering
        \includegraphics[width=0.845\textwidth]{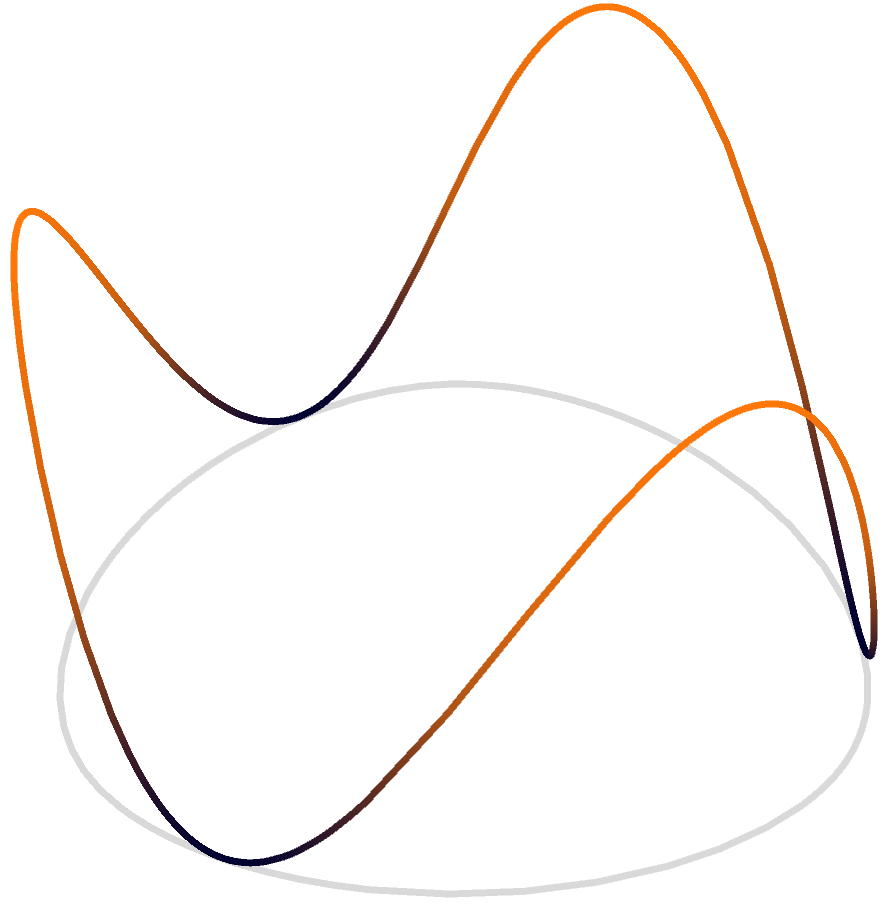} 
        \label{fig:skewloop}
        \caption{A skew loop in $\R^3$, \ \ \ \ \ \ \ \ \ \ with no parallel tangent lines}
    \end{minipage}\hfill
    \begin{minipage}{0.45\textwidth}
        \centering
        \includegraphics[width=0.71\textwidth]{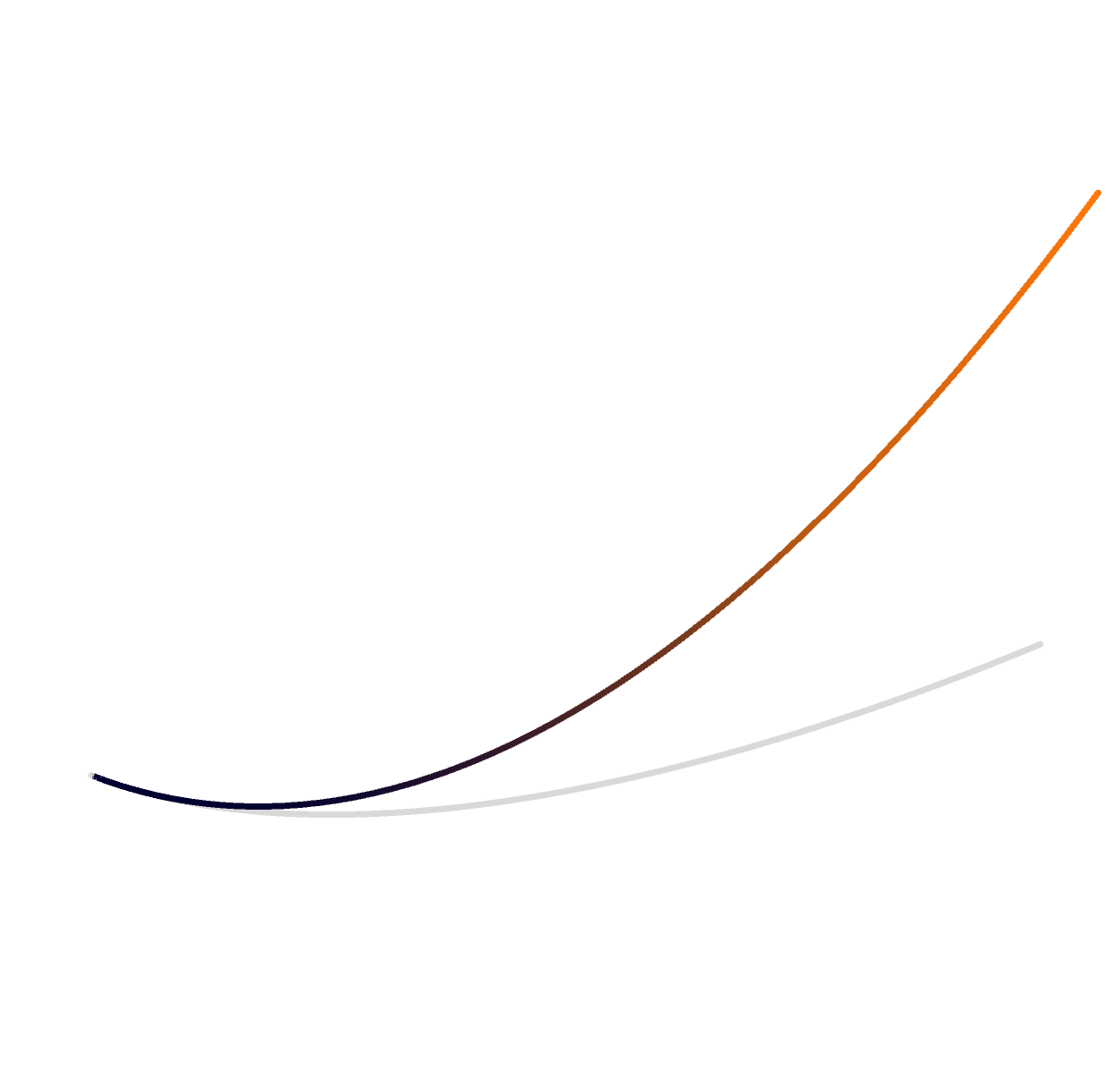} 
        \caption{A totally skew curve in $\R^3$, with no parallel nor intersecting tangent lines}
        \label{fig:skewcurve}
    \end{minipage} \ \ \ \ \ \ \
\end{figure}

Here we consider the following generalization of skew loops.  Throughout this paper, $M$ represents a smooth manifold of dimension $n$.

\begin{definition}  An immersion $f:M \to \R^q$ is \emph{totally nonparallel} if for each pair of distinct points $x,y \in M$, the tangent spaces at $f(x)$ and $f(y)$, considered as affine subspaces of $\R^q$, do not contain parallel lines.  Equivalently, a map is totally nonparallel if each pair of distinct tangent spaces, considered as a pair of linear subspaces of $\R^q$, spans $2n$ dimensions.
\end{definition}

\begin{example}  To supplement the above example of a totally nonparallel immersion $S^1 \to \R^3$, we offer additional examples and one non-example.
\begin{enumerate}
\setlength\itemsep{-.3em}
\item Any map $\R \to \R^2 : x \mapsto (x,g(x))$ with nonvanishing curvature, e.g.\ $x \mapsto (x,x^2)$, is totally nonparallel.
\item The map $\R^2 \to \R^4 : (x,y) \mapsto (x,y,x^2-y^2,2xy)$, induced by $\C \to \C^2: z \mapsto (z,z^2)$, is totally nonparallel.
\item These examples do not generalize to the nonsymmetric quaternionic multiplication: the map $\R^4 \to \R^8$ induced by the quaternionic map $q \mapsto (q,q^2)$ is not totally nonparallel.
\item In Theorem \ref{thm:snbmap} we will see that the first two examples generalize: the graph of any symmetric, bilinear, \emph{nonsingular} map is totally nonparallel.
\end{enumerate}
\end{example}

We say that a map $f : M \to \R^q$ is a $k$-\emph{fold immersion} if, for every $k$-tuple of distinct points $(x_1,\dots,x_k)$ of $M$, the tangent spaces at $f(x_1),\dots,f(x_k)$ span the maximum number of dimensions $kn$.  An example is the moment curve $\R \to \R^k : x \mapsto (x,x^2,\dots,x^k)$, for which the condition is equivalent to the nonvanishing of the determinant of the Vandermonde matrix with distinct entries $(x_1,\dots,x_k)$.  Note that a $1$-fold immersion is an ordinary immersion, and a $2$-fold immersion is a totally nonparallel immersion.  

Our study of totally nonparallel immersions, and the more general $k$-fold immersions, follows a recent surge in the study of certain maps $M \to \R^q$; specifically, those which are distinguished by conditions on the \emph{configuration space} $F_k(M) \coloneqq \left\{ (x_1,\dots,x_k) \ \big| \ x_i \neq x_j \mbox{ for } i \neq j \right\}$ of $k$-tuples of distinct points of $M$.  For example, a smooth map $f: M \to \R^q$ is called:

\begin{itemize}
\item \emph{$k$-regular} if for each $x \in F_k(M)$, the collection of vectors $\left\{f(x_1), \dots, f(x_k) \right\}$ is linearly independent (\cite{BlagojevicLueckZiegler}, \cite{BoltjanskiiEtal}, \cite{Borsuk}, \cite{BuczynskiEtal}, \cite{Chisholm},
\cite{CohenHandel}, \cite{Handel}, \cite{Handel2}, \cite{Handel3},
\cite{HandelSegal}).
\item \emph{$k$-skew} if for each $x \in F_k(M)$, the collection of affine subspaces $\left\{T_{f(x_1)}\R^q, \dots, T_{f(x_k)}\R^q\right\}$ of $\R^q$ is affinely independent (\cite{Baralic}, \cite{BaralicEtal}, \cite{BlagojevicLueckZiegler}, \cite{GhomiTabachnikov}, \cite{Stojanovic}).
\item \emph{$k$-neighborly} if for each $x \in F_k(M)$, there is an affine hyperplane $H$ supporting $f(M)$ which touches $f(M)$ at exactly the points $f(x_1), \dots, f(x_k)$; that is, $H$ contains the points $f(x_1), \dots, f(x_k)$ and one of the open half-spaces determined by $H$ contains the remaining points of $f(M)$ (\cite{KalaiWigderson}, \cite{Vassiliev}).
\end{itemize}
Complex $k$-regular embeddings \cite{BlagojevicEtal}, $\ell$-skew-$k$-regular embeddings (\cite{BlagojevicLueckZiegler}, \cite{Stojanovic}), tangent-bundle embeddings (\cite{Ghomi}, \cite{StojanovicTabachnikov}), and skew branes \cite{TabachnikovBranes}, \cite{TabachnikovBranes2}, \cite{ShaSolomon}) have also been studied recently.

The main focus of our paper is the following question.
\begin{mainquestion}
Given a smooth manifold $M$, what is the minimum dimension $\TN(M)$ such that there exists a totally nonparallel immersion $M \to \R^{\TN(M)}$?
\end{mainquestion}

The Main Question has been studied for all of the conditions listed above, though the majority of results focus only on lower bounds.  In general, lower bounds may be computed as follows.  A map $M \to \R^q$ satisfying one of the above $k$-tuple conditions allows for the construction of a corresponding $S_k$-equivariant map, either on the configuration space itself or on a vector bundle over it, satisfying a certain nondegeneracy condition.  For example, a $k$-regular embedding $M \to \R^q$ induces an $S_k$-equivariant map from $F_k(M)$ to the Stiefel manifold $V_k(\R^q)$, and we will see that a $k$-fold immersion $M \to \R^q$ induces an immersion of the unordered configuration space $F_k(M) / S_k \to \R^q$.  In either of these cases, and also for the other $k$-tuple conditions above, obstructions may be computed by Stiefel-Whitney class calculations or other methods.  These obstructions are often fairly difficult to compute, and for the majority of the conditions on the above list, are computed only in the case $M = \R^n$.

A notable exception is the $2$-skew embedding condition, introduced and studied as ``totally skew embeddings'' by Ghomi and Tabachnikov \cite{GhomiTabachnikov}.  An embedding $f : M \to \R^q$ is called \emph{totally skew} if for each pair of distinct points $x$, $y$ of $M$, the tangent spaces at $f(x)$ and $f(y)$ neither intersect nor contain parallel directions.  The simplest example is the \emph{Veronese embedding} of $\R \to \R^3$ given by $x \mapsto (x,x^2,x^3)$, depicted above in Figure \ref{fig:skewcurve}.  Ghomi and Tabachnikov sought the minimum dimension $\TS(M)$ such that there exists a totally skew embedding $M \to \R^{\TS(M)}$ and found that obstructions in the case $M = \R^n$ are related to the generalized vector field problem and to immersions of real projective spaces.

The search for $\TS(M)$ was continued by Barali\'{c}, et al.\ \cite{BaralicEtal}, who employ Stiefel-Whitney class calculations to compute lower bounds for $\TS(M)$ when $M$ is a projective space, a Grassmann manifold, or (separately by Barali\'{c} \cite{Baralic}) a certain quasitoric manifold.  The best estimates arise in the case of real projective spaces in dimensions which are powers of $2$: if $n = 2^m$, then $4n - 1 \leq \TS(\RP^n) \leq 4n+1$; the upper bound is by a genericity argument given by Ghomi and Tabachnikov.

Due to the similarities between the conditions, all of the above arguments may be reproduced for totally nonparallel immersions to obtain similar bounds; for example, the argument of \cite{BaralicEtal} and a genericity argument yield $4n-2 \leq \TN(\RP^n) \leq 4n$ when $n$ is a power of $2$.  However, the authors of \cite{BaralicEtal} do not take advantage of a certain $\Z_2$-equivariance present in the obstruction, and so these lower bounds for $\TS$ and $\TN$ are not sharp.  Specifically, in the totally nonparallel case, the obstruction takes the following form.

\begin{observation}
\label{obs:main} A totally nonparallel immersion $f : M \to \R^q$ induces an immersion of the unordered configuration space $F_2(M) / \Z_2 \to \R^q$, given by $[(x,y)] \mapsto f(x) + f(y)$.
\end{observation}

Following a discussion of this problem with Don Davis, he was able to show the following non-immersion theorem, from which follows the improved lower bound $4n - 1 \leq \TN(\RP^n)$ when $n = 2^m$.

\begin{thm}[Davis \cite{Davis2019}] \label{thm:davis} For $n = 2^m$, the unordered configuration space $F_2(\RP^n) / \Z_2$ does not immerse in $\R^{4n-2}$.
\end{thm}

The theorem of Davis also applies to totally skew embeddings, yielding the improved lower bound $4n \leq \TS(\RP^n)$ for $n = 2^m$.

The goal of this paper is the exact value of $\TN(\RP^n)$.

\begin{thm}
\label{thm:tnrpn}
For $n = 2^m$, $\TN(\RP^n) = 4n-1$.
\end{thm}

This theorem provides the first answer to the Main Question for any closed manifold $M \neq S^1$ for any of the conditions listed above.  In particular, the only existence results for any of the conditions above arise from generic bounds or constructions which are not likely to be sharp.  One exception is for $k$-regular embeddings: strong existence results were computed using methods of algebraic geometry \cite{BuczynskiEtal}, though still these results only hold for $M = \R^n$.  We do not know if the methods there could improve our upper bounds for $\TN(\R^n)$, which are listed in Table \ref{tab:tcimmersions} in Section \ref{sec:localex}.

We attack the existence problem in this paper using a combination of singularity theory and h-principle techniques, which culminate in a Whitney-type theorem for totally nonparallel immersions.  Recall that the weak Whitney theorems state that a generic map $f: M \to \R^{2n+1}$ is an embedding, and that a generic map $f: M \to \R^{2n}$ is an immersion, hence a local embedding.  Then the Whitney trick may be used to prove the strong Whitney embedding theorem: that there exists an embedding $M \to \R^{2n}$.  Here we will employ a similar strategy.  We will first see that a generic map $f : M \to \R^{4n}$ is totally nonparallel.  Then, after a detour to understand the space of symmetric bilinear \emph{nonsingular} maps, we show that a generic map $f : M \to \R^{4n-1}$ is locally totally nonparallel.  Finally, we use the removal of singularities h-principle technique, due to Gromov and Eliashberg \cite{GromovEliashberg}, to systematically modify a locally totally nonparallel immersion until it is (globally) totally nonparallel.

Theorem \ref{thm:tnrpn} then follows by combining the theorem of Davis with the following: 

\begin{thm} 
\label{thm:main}
Every smooth $n$-dimensional manifold $M$ admits a totally nonparallel immersion into $\R^{4n-1}$.
\end{thm}

The removal-of-singularities technique has been used in similar settings, but there has not appeared a coherent exposition of the technical details.  In particular, the idea above matches the general structure of \cite{szucs}, in which Sz\H{u}cs used the Gromov-Eliashberg technique to offer an alternate proof of the $N = \R^q$ case of Haefliger's Theorem on embeddings (stated below).  In the present paper we address and remedy certain technical gaps from \cite{szucs}.

\begin{thm}[Haefliger \cite{Haefliger2}]  Let $M$ be a smooth closed manifold of dimension $n$ and $N$ a smooth manifold of dimension $q$, where $q > \frac32 n + \frac32$.  Then the existence of a differentiable embedding $M \to N$ is equivalent to the existence of a $\Z_2$-equivariant map $g : M \times M \to N \times N$ such that $g^{-1}(\Delta N) = \Delta M$.
\end{thm}

We believe that a Haefliger-type h-principle statement holds for totally nonparallel immersions.  Our proposed statement is the converse of  Observation \ref{obs:main}, assuming the appropriate range of dimensions:

\begin{conjecture}
\label{con:main}
Let $M$ be a smooth closed manifold of dimension $n$ and let $q \geq \frac72 n$.  If there exists an immersion of the unordered configuration space $g : F_2(M) / \Z_2 \to \R^q$, then there exists a totally nonparallel immersion $M \to \R^q$.
\end{conjecture}

In the case $q = 4n-1$, the existence of an immersion is not an additional restriction, since every $2n$-manifold immerses in $\R^{4n-1}$.  Nevertheless, our proof of Theorem \ref{thm:main} follows the flavor of this statement, in the sense that we begin with a certain immersion $F_2(M) / \Z_2 \to \R^{4n-1}$ and manipulate it to build a totally nonparallel immersion.  The majority of our construction still applies with the more general bound $q \geq \frac72 n$.  In Section \ref{sec:conclusion} we discuss which ingredients are missing from a complete proof.

More generally, it is our hope to develop a framework which extends the results of this paper and applies to general configuration space conditions, especially those listed above.  This is a serious endeavor; the suggestion to develop h-principle results for such configuration space conditions appeared over thirty years ago in the book of Gromov \cite{Gromov} (pages 51--52, Question and subsequent examples), yet there has still been no progress in this area.

The paper is organized as follows.  We begin with results for totally nonparallel immersions which are similar in flavor to those obtained by Ghomi and Tabachnikov \cite{GhomiTabachnikov} for totally skew embeddings.  In particular, Sections \ref{sec:localobs}--\ref{sec:localex} establish relationships among totally nonparallel immersions, the generalized vector field problem, the immersion and embedding problems for real projective spaces, and symmetric bilinear nonsingular maps.  We continue in Sections \ref{sec:review}--\ref{sec:removal} with a brief review of standard notions in singularity theory which will be used throughout the paper, as well as an illustrative demonstration of the removal-of-singularities h-principle technique.  In Section \ref{sec:removaltn} we offer a very rough outline for the proof of Theorem \ref{thm:main}, which will pave the way for the technicalities addressed in future sections.

We continue in Section \ref{sec:local} with a treatment of the local totally nonparallel problem, by revisiting the space of symmetric bilinear nonsingular maps and studying the structure of this space in detail.  In Section \ref{sec:outline}, we refine the outline given in Section \ref{sec:removaltn}, incorporating the local considerations of Section \ref{sec:local}, and we prove the theorem under certain technical assumptions which are verified in Section \ref{sec:global}.  This verification process hinges on certain transversality theorems, designed specially for our needs, which are proven in Section \ref{sec:transversality}.  We conclude with some discussion on various possible improvements of the results of the paper. \\

\textbf{Acknowledgments}  I am grateful for many useful discussions during the preparation of this article, especially to: Martin Bendersky, Pavle Blagojevi\'{c}, Donald Davis, Yasha Eliashberg, Andr\'{a}s Juh\'{a}sz, Kee Lam, Daniel Levine, and Sergei Tabachnikov.  I would also like to thank Peter Landweber and an anonymous referee for useful comments on the presentation.

\subsection{Obstructions to the existence of totally nonparallel immersions}
\label{sec:localobs}

We begin with simple lower bounds for $k$-fold immersions, and proceed to stronger bounds for totally nonparallel immersions, establishing strong relationships among the notions of totally nonparallel immersions, the generalized vector field problem, and the immersion and embedding problems for real projective spaces.  Many of the results of this section are similar to results obtained by Ghomi and Tabachnikov for totally skew embeddings \cite{GhomiTabachnikov}, and we do not give proofs for those statements which follow from theirs with only superficial adaptations.

If there exists a $k$-fold immersion $M \to \R^q$, then $q \geq kn$, since by definition, the tangent spaces at each $k$ distinct points collectively span a space of dimension $kn$ in the target space.  For compact manifolds $M$, this bound may be slightly improved with a Morse theory argument.

\begin{prop}
Let $M$ be a smooth compact $n$-dimensional manifold.  If $M$ admits a totally nonparallel immersion into $\R^q$, then $q \geq 2n+1$.
\end{prop}

\begin{proof}
Choose a generic height function in $\R^q$.  The restriction to $f(M)$ has (at least) two critical points, for which the tangent spaces must span a $(2n)$-dimensional space orthogonal to the chosen height direction.
\end{proof}

This idea generalizes to $k$-fold immersions for compact manifolds $M$ whose Betti number sum is at least $k$.  In any case, stronger lower bounds are given as follows.

\begin{prop}
\label{prop:skimm}
Let $M$ be a smooth $n$-dimensional manifold.  If $M$ admits a $k$-fold immersion into $\R^q$, then there is an immersion of the unordered configuration space $F_k(M) / S_k$ into $\R^q$.
\end{prop}

\begin{proof}
Let $f : M \to \R^q$ be a $k$-fold immersion.  Then the map $g : F_k(M) \to \R^q$, which maps $x = (x_1,\dots,x_k)$ to the sum $f(x_1) + \cdots + f(x_k)$, is $S_k$-equivariant with respect to the trivial action on $\R^q$.  Moreover, $g$ is an immersion, since the tangent space at $g(x)$ is the direct sum of the tangent spaces at the $k$ image points $f(x_1), \dots, f(x_k)$, and hence has dimension $kn$.  Thus $g$ induces an immersion of the unordered configuration space $F_k(M) / S_k$ to $\R^q$.
\end{proof}

\begin{remark}
Using the proposition above, obstructions to $k$-fold immersions may be computed in terms of Stiefel-Whitney classes of unordered configuration spaces; this is the content of Davis' proof of Theorem \ref{thm:davis}.  For $k$-regular, $\ell$-skew, and $k$-regular-$\ell$-skew embeddings, similar computations were done by Blagojevi{\'c}, L{\"u}ck, and Ziegler \cite{BlagojevicLueckZiegler}; however, their computations use a result of Hung \cite{Hung} on the mod-$2$ equivariant cohomology algebras of configuration spaces, which contains an error.  The error was recently reconciled by Blagojevi{\'c}, Cohen, Crabb, L{\"u}ck, and Ziegler, and the updated lower bounds may be found in \cite{BlagojevicCohenEtal}, Theorems 5.14, 5.18, and 5.22.
\end{remark}

In the case of totally nonparallel immersions, we also have the following more specific results, which connect the theory of totally nonparallel immersions with the generalized vector field problem and the immersion problem for real projective spaces.  These results are similar to those obtained by Ghomi and Tabachnikov in their study of totally skew embeddings, and the proofs are nearly identical, so we do not repeat them here.  It is worth mentioning, however, that although the totally skew condition is strictly stronger than the totally nonparallel condition, one does not obtain a stronger version of Proposition \ref{prop:immproj} for totally skew embeddings.

\begin{prop}
\label{prop:generalizedvector}
If $\R^n$ admits a totally nonparallel immersion into $\R^q$, then there exist $n$ linearly independent sections of the bundle $(q-n)\xi_{n-1}$, where $\xi_{n-1} \to \RP^{n-1}$ is the tautological line bundle.
\end{prop}

\begin{cor}
For $n \neq 1, 2,4,8$, there is no totally nonparallel immersion $\R^n \to \R^{2n}$.
\end{cor}

A compelling open question is whether there exist totally nonparallel immersions $\R^4 \to \R^8$ or $\R^8 \to \R^{16}$.

\begin{prop}
\label{prop:immproj}
If $\R^n$ admits a totally nonparallel immersion into $\R^q$ and $n \neq 1,2,4,8$, then there exists an immersion of $\RP^{n-1}$ into $\R^{q-n-1}$.
\end{prop}

The lower bounds resulting from the above statements are summarized in Table \ref{tab:tcimmersions}, following a discussion of the existence problem.

\subsection{Existence of totally nonparallel and $k$-fold immersions of $\R^n$}
\label{sec:localex}

We begin our study of the existence problem with a genericity statement for $k$-fold immersions.  Here, a subset $X \subset Y$ is \emph{residual} if it is a countable intersection of open dense subsets of $Y$.

\begin{thm}
\label{thm:gen}
For $q \geq 2kn$, the subset $\left\{ f \in C^\infty(M,\R^q) \ \big| \ f \mbox{ is a } k\mbox{-fold immersion}\right\}$ is a residual subset of $C^\infty(M,\R^q)$ in the Whitney $C^\infty$ topology.
\end{thm}

Intuitively, the $k$-fold immersion condition corresponds to the configuration space condition that the map $F_k(M) / S_k \to \R^q : (x_1, \dots, x_k) \mapsto f(x_1) + \cdots + f(x_k)$ is an immersion, which we can expect generically for $q \geq 2kn = 2\dim(F_k(M) / S_k)$.  We formalize this in Section \ref{sec:review}, following the statement of the Multijet Transversality Theorem (Theorem \ref{thm:mtt}).

Focusing now only on totally nonparallel immersions, we find that stronger upper bounds for the number $\TN(\R^n)$ arise by construction via nonsingular symmetric bilinear maps.  Here, a symmetric bilinear map $B : \R^n \times \R^n \to \R^{q-n}$ is \emph{nonsingular} if $B(x,y) = 0$ only when $x=0$ or $y=0$.  In that sense, totally nonparallel immersions are related not only to the immersions of real projective space as in Proposition \ref{prop:immproj}, but also to embeddings of real projective space; since nonsingular symmetric bilinear maps are known to produce such embeddings (see \cite{James4}).

\begin{thm}
\label{thm:snbmap}
Let $B: \R^n \times \R^n \to \R^{q-n}$ be a symmetric, bilinear map, and let $Q : \R^n \to \R^{q-n} : x \mapsto B(x,x)$ be the associated quadratic map.  Then the graph of $Q$ is a totally nonparallel immersion if and only if $B$ is nonsingular.
\end{thm}

The backwards direction will follow from the more general Proposition \ref{prop:localtotallynonparallel}.  The proof there can be slightly modified, so that only quadratic maps $(x,B(x,x))$ are considered, to establish the forward direction.

\begin{cor}
\label{cor:quad}
There exists a totally nonparallel immersion of $\R^n$ to $\R^{3n-1}$, and when $n$ is even, of $\R^n$ to $\R^{3n-2}$.
\end{cor}

\begin{proof}
Define the quadratic map
$$Q : \R^n \to \R^{2n-1} : (x_0,\dots,x_{n-1}) \mapsto \bigg( \sum_{i+j =0} x_ix_j , \sum_{i+j =1} x_ix_j , \dots, \sum_{i+j = 2n-2} x_ix_j \bigg).$$ 
The associated symmetric bilinear map $\R^n \times \R^n \to \R^{2n-1}$ corresponds to multiplying two degree-$(n-1)$ polynomials, a nonsingular operation.  When $n$ is even, one may consider the similarly defined complex map.
\end{proof}

\begin{cor}
\label{cor:noquad}
There is no quadratic totally nonparallel immersion $\R^3 \to \R^7$.
\end{cor}

\begin{proof}
A quadratic totally nonparallel immersion $\R^3 \to \R^7$ would correspond to a nonsingular symmetric bilinear map $\R^3 \times \R^3 \to \R^4$, which would induce an embedding $\RP^2 \to \R^3$, but no such embedding exists.
\end{proof}

For all $n > 2$, there are still gaps between the lower bounds for $\TN(\R^n)$ obtained by Proposition \ref{prop:generalizedvector} and the upper bounds for $\TN(\R^n)$ of Corollary \ref{cor:quad}, as indicated in the table below.  We obtained the lower bounds from Proposition \ref{prop:generalizedvector} by consulting Table 5.14 of \cite{LamRandall}.

\begin{align*}
{\setlength{\arraycolsep}{.38em}
\begin{array}{c||c|c|c|c|c|c|c|c|c|c|c|c|c|c|c|c|c}
n \ & 1 & 2 & 3 & 4 & 5 & 6 & 7 & 8 & 9 & 10 & 11 & 12 & 13 & 14 & 15 & 16 & 17 \\
\hline \hline
\TN(\R^n) \geq & \bf{2} & \bf{4} & 7 & 8 & 13 & 14 & 15 & 16 & 25 & 26 & 28 & 29 & 32 & 37 & 38 & 39 & 49 \\ 
\TN(\R^n) \leq \ & \bf{2} & \bf{4} & 8 & 10 & 14 & 16 & 20 & 22 & 26 & 28 & 32 & 34 & 38 & 40 & 44 & 46 & 50 \\
\end{array}}
\end{align*}
\begin{table}[!h]
\vskip-.7cm
\caption{Lower and upper bounds for $\TN(\R^n)$}
\label{tab:tcimmersions}
\end{table}

We take a moment to call attention to the ``jumps'' which occur in the lower bounds whenever $n - 1$ is a power of $2$.

In the spirit of Corollary \ref{cor:noquad}, if it were known that any totally nonparallel immersion produced a \emph{quadratic} totally nonparallel immersion (say, by considering only the $2$-jet), then totally nonparallel immersions would produce \emph{embeddings} of real projective spaces, and the bounds resulting from Proposition \ref{prop:immproj} could be significantly improved.  In the conclusion of this paper, we discuss why we cannot necessarily expect this phenomenon to occur.

\subsection{Review of singularity theory and relevant definitions}
\label{sec:review}

Here we review the main tools and introduce notation which we will use in the remainder of this paper.  We refer to Chapters II and VI of the book by Golubitsky and Guillemin \cite{GolubitskyGuillemin} for the bulk of the results below.

Let $X$ and $Y$ be smooth manifolds, with respective dimensions $n$ and $q$, and let $q > n$.  We assume familiarity with the $s$-jet space $J^s(X,Y)$.  For a comprehensive introduction to jet bundles, we recommend \cite{Saunders}.

\begin{thm}[Thom Transversality Theorem (see \cite{GolubitskyGuillemin} Theorem II.4.9)]  Let $X$ and $Y$ be smooth manifolds and $W$ a stratified subset of $J^s(X,Y)$.  Then the set
\[
\left\{ f \in C^\infty(X,Y) \ | \ j^sf \pitchfork W \right\}
\]
is a residual subset of $C^\infty(X,Y)$  in the Whitney $C^\infty$ topology.
\end{thm}

In practice, this theorem is quite often applied to sets of the following form.  Consider the space $\Hom(\R^n,\R^q)$ of linear maps from $\R^n$ to $\R^q$, where $q > n$.  Let $L^r(\R^n,\R^q)$ denote the subset of $\Hom(\R^n,\R^q)$ consisting of those linear maps with corank $r$; that is, rank $n - r$.  Then $L^r(\R^n,\R^q)$ is a submanifold of $\Hom(\R^n,\R^q)$ with codimension $(q-n+r)r$.

The $1$-jet bundle $J^1(X,Y)$ has the local trivialization $\R^n \times \R^q \times \Hom(\R^n,\R^q)$.  Define the subfiber-bundle $S_r \subset J^1(X,Y)$, which has fiber $L^r(\R^n,\R^q)$ over each point $(x,z)$.  Then $S_r$ has codimension $(q-n+r)r$.

A smooth map $f : X \to Y$ is an immersion if its differential $df_x : T_xX \to T_zY$  has rank $n$ everywhere.  This occurs if and only if the $1$-jet extension $j^1f$ lies in the open submanifold $S_0$, so the space of immersions is open by the definition of $C^\infty$ topology.  On the other hand, $f$ is an immersion if and only if the $1$-jet extension $j^1f$ is disjoint from $S_r$ for every $r>0$.  By the Thom Transversality Theorem, the subset $\left\{ f \in C^\infty(X,Y) \ | \ j^1f \pitchfork S_r \right\}$ is dense.  When $q \geq 2n$, $\codim(S_r) = (q-n+r)r > n$, and transversality $j^1f \pitchfork S_r$ is disjointness.  Hence when $q \geq 2n$, the space of immersions from $X \to Y$ is dense in $C^\infty(X,Y)$.  

Now consider the $k$-fold product $X^k$ and the configuration space $F_k(X)$ of $k$-tuples of distinct points of $X$.  Let $\alpha : J^s(X,Y) \to X$ be the projection, and let $\alpha^k : J^s(X,Y)^k \to X^k$ be the induced projection.  We define the $k$-fold $s$-multijet bundle $J^s_k(X,Y) = (\alpha^k)^{-1}(F_k(X))$, and we define the $k$-fold $s$-multijet of $f : X \to Y$ by $j_k^sf : F_k(X) \to J^s_k(X,Y) : (x_1,\dots,x_k) \mapsto (j^sf(x_1),\dots,j^sf(x_k))$. 

\begin{thm}[Multijet Transversality Theorem (see \cite{GolubitskyGuillemin} Theorem II.4.13)] Let $X$ and $Y$ be smooth manifolds with $W$ a submanifold of $J_k^s(X,Y)$.  Then
\[
\left\{ f \in C^\infty(X,Y) \ | \ j_k^s f \pitchfork W \right\}
\]
is a residual subset of $C^\infty(X,Y)$  in the Whitney $C^\infty$ topology.
\label{thm:mtt}
\end{thm}

As an application of the multijet transversality theorem, we give the proof of Theorem \ref{thm:gen}.

\begin{proof}[Proof of Theorem \ref{thm:gen}]
For any smooth map $f: M \to \R^q$ and $k$-tuple of distinct points $(x_1,\dots,x_k)\in F_k(M)$, we may write the corresponding $1$-multijet extension of $f$ by
\[
j^1_k f(x_1,\dots,x_k) \coloneqq (x_1,\dots,x_k,f(x_1),\dots,f(x_k),df_{x_1},\dots,df_{x_k}).
\]
We define the singularity set $S$ in the $k$-fold $1$-multijet $J_1^k$ as follows:
\[
S \coloneqq \left\{ (x_1,\dots,x_k,y_1,\dots,y_k,P_1,\dots,P_k) \in J^1_k \ \big| \ \operatorname{dim}\left(\operatorname{Span} \left\{\operatorname{Im}(P_1),\dots,\operatorname{Im}(P_k)\right\}\right) < kn \right\}.
\]
A map $f: M \to \R^q$ is a $k$-fold immersion if and only if the $1$-multijet extension of $f$ is disjoint from $S$.  We would like to show that $S$ is a stratified subset of $J_k^1(M,\R^q)$ of codimension $q-kn+1$, for then it follows from the Multijet Transversality Theorem that the $1$-multijet extension of a generic $f \in C^\infty(M,\R^q)$ intersects $S$ transversely, hence is disjoint from $S$ when $q \geq 2kn$.

The bundle $J^1_k \to F_k(M)$ has the local trivialization $(\R^n)^k \times (\R^q)^k \times (\Hom(\R^n,\R^q))^k$.  Since $S$ is independent of the source $x$ and the target $y$, it is enough to compute the codimension of $S_{(x,y)} \subset (J^1_k)_{(x,y)} \simeq (\Hom(\R^n,\R^q))^k$.  For this, consider the projection $\pi_i$ of $(\R^n)^k$ onto the $i$th factor of $\R^n$, and observe that there is a diffeomorphism
\[
\Psi : (\Hom(\R^n,\R^q))^k \to \Hom((\R^n)^k,\R^q) = \Hom(\R^{kn},\R^q),
\]
where $\Psi(P_1,\dots,P_k)$ is defined by its action on $v \in \R^{kn}$ as follows:
\[
\Psi(P_1,\dots,P_k)(v) = \sum_{i=1}^k P_i(\pi_i(v)).
\]
Let $L^r \subset \Hom(\R^{kn},\R^q)$ be the subset of linear maps which have corank $r$.  Then for any source $x$ and target $y$, $(P_1,\dots,P_k) \in S_{(x,y)}$ if and only if $\Psi(P_1,\dots,P_k) \in L^r$ for some $r > 0$.  By the corank formula above, $\operatorname{codim}(L^r) = (q-kn+r)r$, so the union of the $L^r$ is a stratified set with codimension $q-kn+1$.  It follows that $S \subset J^1_k$ is stratified with the same codimension, and the proof is complete.
\end{proof}

\subsection{The removal of singularities h-principle technique}
\label{sec:removal}

We offer a brief example to illustrate the ``removal of singularities'' h-principle technique due to Gromov and Eliashberg \cite{GromovEliashberg}.  Following the idea of the treatment in Gromov's book \cite{Gromov}, page 49, we offer a proof of the following special case of the Smale-Hirsch immersion theorem.  

\begin{thm}[Smale-Hirsch Theorem, special case]
\label{thm:smalehirsch}
Let $M$ be a smooth closed manifold of dimension $n$, and let $q > n$.  If there exists a bundle monomorphism $TM \to \R^q$, then there exists an immersion $M \to \R^q$.
\end{thm}

\begin{proof}
A bundle monomorphism $TM \to \R^q$ may be considered as $q$ $1$-forms $(\omega_1,\dots,\omega_q)$ which span the cotangent bundle.  If the $1$-forms were exact; that is, if for all $i$, $\omega_i = df_i$ for some function $f_i : M \to \R$, then we would have the desired immersion $f = (f_1,\dots,f_q)$.  The strategy is to replace the $1$-forms, one-by-one, with exact forms, without disrupting the feature that the $1$-forms span the cotangent bundle.

Let us first consider $\omega_1$, and suppose that naively we replace $\omega_1$ with $df_1$ for some arbitrary smooth function $f_1 : M \to \R$.  This replacement would cause issues only on the set
\[
\Sigma_1 \coloneqq \left\{ x \in M \ \big| \ \rank((\omega_2,\dots,\omega_q)) < n \right\},
\]
since on the complement of $\Sigma_1$, $(df_1,\omega_2,\dots,\omega_n)$ has full rank regardless of how $f_1$ is defined.  Thus we may perform the replacement by defining $f_1$ carefully on $\Sigma_1$, and then we may extend arbitrarily to $\Sigma_1^c$.  We use the notation $\omega = (\omega_1,\dots,\omega_q)$ and $\tilde{\omega} = (\omega_2,\dots,\omega_q)$.

Since $\omega$ has rank $n$ everywhere, $\tilde{\omega}$ has rank at least $n-1$ everywhere, so $\ker(\tilde{\omega})$ is a $1$-dimensional line bundle on $\Sigma_1$.  Moreover, the line bundle is trivial, since $\omega_1(\ker(\tilde{\omega}))$ is nonzero on $\Sigma_1$.

We consider the special case in which $\ker(\tilde{\omega})$ is transverse to $\Sigma_1$.  In this case, we may consider the trivial line bundle embedded in $M$ as $\Sigma_1 \times [-\varepsilon,\varepsilon]$.  We may define $f_1 = 0$ on $\Sigma_1$ and extend $f_1$ in the direction of the line bundle in such a way that $df_1(\ker(\tilde{\omega}))  =\omega_1(\ker(\tilde{\omega}))$, and then arbitrarily extend $f_1$ to all of $M$ (see Figure \ref{fig:sigma}).

\begin{figure}[h!t]
    \centering
        \includegraphics[width=0.9\textwidth]{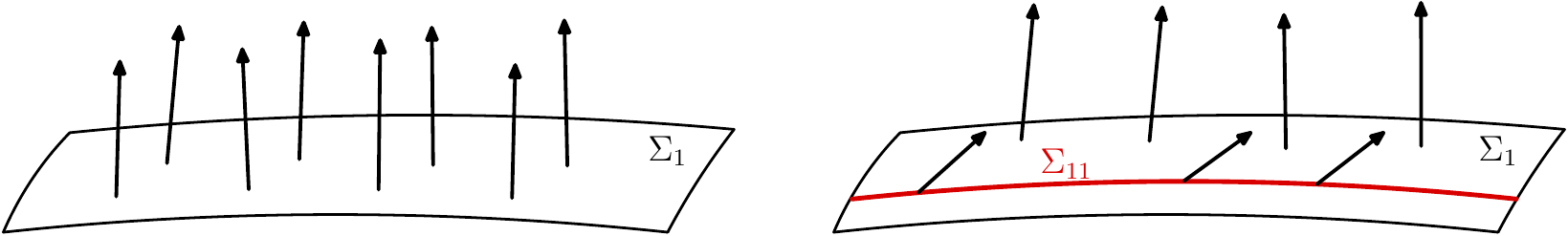} 
        \caption{Two depictions of the singularity set $\Sigma_1$ with a section of $\ker(\tilde{\omega})$.  Left: The special case in which $\ker(\tilde{\omega})$ is transverse to $\Sigma_1$.  Right: The (less) special case, in which $\ker(\tilde{\omega})$ is transverse to $\Sigma_{11}$.}
        \label{fig:sigma}
\end{figure}

If we are not in the special case, we consider the set where the transversality to $\Sigma_1$ fails:
\[
\Sigma_{11} \coloneqq \left\{ x \in \Sigma_1 \ \big| \ \ker(\tilde{\omega}) \mbox{ is tangent to } T_x\Sigma_1 \right\},
\]
and we consider the special case in which the restriction of the kernel to $\Sigma_{11}$ is transverse to $\Sigma_{11}$.  In this case, we may define $f_1$ first on $\Sigma_{11}$, then use the transversality to extend to a neighborhood of $\Sigma_{11}$ in $\Sigma_1$, such that $df_1(\ker(\tilde{\omega}))  =\omega_1(\ker(\tilde{\omega}))$ for points of $\Sigma_{11}$, then extend the function $f$ arbitrarily to all of $\Sigma_1$, then use transversality away from $\Sigma_{11}$ to extend $f_1$ to a neighborhood of $\Sigma_1$, and finally extend $f_1$ to all of $M$.

The spaces $\Sigma_{1\dots 1}$ can be defined similarly; these are \emph{Thom-Boardman singularity spaces}, and the Thom-Boardman Theorem guarantees that for generic $\tilde{\omega}$, these provide a stratification of $\Sigma_1$.  Thus we can always define $f_1$ in a manner resembling the above procedure: we define $f_1$ on the lowest-dimensional stratum and inductively from there; beginning, if necessary, with a small perturbation of $\omega$ within the space of full rank maps, to guarantee genericity.  After $f_1$ is replaced, we continue in a similar manner, replacing each $\omega_i$ as above, using perturbations of the remaining $\omega_j$, $j > i$, and of the already-replaced $f_j$, $j < i$, whenever necessary, until all $q$ $1$-forms have been replaced.
\end{proof}

In Section \ref{sec:step2}, the Thom-Boardman spaces will make a second appearance, and we will discuss the \emph{intrinsic derivative} method of Porteous for dealing with these singularity spaces.

\subsection{Removal of singularities for totally nonparallel immersions}
\label{sec:removaltn}

In this final introductory section, we discuss the broad idea for the proof of Theorem \ref{thm:main} and the difficulties encountered in using the removal of singularities technique.  In the proof above, at each stage of the replacement we must make a perturbation of $\omega$ to guarantee a certain generic property.  This perturbation is safe because the full rank condition is open.  When attempting to repeat the argument for totally nonparallel immersions, a small perturbation could disrupt important properties near the diagonal of $M \times M$, and so appropriate modifications must be made.  We follow the idea of \cite{szucs}, which contains a similar argument for Haefliger's embedding theorem.

\begin{notation} Given a map $f : M \to \R^q$, we define $f \oplus f : M \times M \to \R^q : (x,y) \mapsto f(x) + f(y)$.
\end{notation}

We recall the statement of Proposition \ref{prop:skimm} in the case $k=2$:  if $f : M \to \R^q$ is a totally nonparallel immersion, then the map $g = f \oplus f : M \times M \to \R^q$ restricts to a $\Z_2$-equivariant immersion of $F_2(M)$ into $\R^q$.

In the spirit of Conjecture \ref{con:main}, suppose conversely that we begin with a $\Z_2$-equivariant map $g : M \times M \to \R^q$ which restricts to an immersion on $F_2(M)$.  Write $g$ in components $g = (g_1,\dots,g_q)$, where each $g_i$ is a real-valued map.  Note that if each $g_i$ happens to be of the form $f_i \oplus f_i$ for some map $f_i : M \to \R$, then $f$ is the desired totally nonparallel immersion.  The strategy is to replace the functions $g_i$, one-by-one, with functions of the form $f_i \oplus f_i$ for some $f_i : M \to \R$.  If we are able to perform these replacements, such that the resulting map is still an immersion on $F_2(M)$, then $f = (f_1,\dots,f_q)$ is a totally nonparallel immersion.  

We start with $g_1$.  Our primary concern during the replacement is that the differential of the new map $(f_1 \oplus f_1, g_2, \dots, g_q)$ must have full rank $2n$ at every point of $F_2(M)$.  This is not an issue at points $(x,y) \in F_2(M)$ at which $\tilde{g} = (g_2,\dots,g_q)$ already has rank $2n$; indeed, at such points we may replace $g_1$ arbitrarily without causing a problem.  However, we must be more careful on the set
\[
\Sigma_{\tilde g} \coloneqq \left\{ (x,y) \in F_2(M) \ \big| \ \rank(d\tilde g_{(x,y)}) = 2n - 1 \right\}.
\]
By definition of $\Sigma_{\tilde g}$, at each point $(x,y)$ there exists a one-dimensional kernel $\ker(d\tilde g_{(x,y)})$, which we consider as a line subbundle of $T(F_2M)\big| \Sigma_{\tilde g}$.  Since $dg$ has full rank on $F_2(M)$, $d(g_1)$ does not vanish in the direction of this kernel, so $\ker(d\tilde g)$ is a trivial line bundle.  In particular, let $s$ be the nonzero section on which $d(g_1)$ is identically $1$.  If we can find $f_1 : M \to \R$ such that $d(f_1 \oplus f_1)_{(x,y)}(s(x,y)) = 1$ for all $(x,y) \in \Sigma_{\tilde{g}}$, then the new map $(f_1 \oplus f_1, g_2, \dots, g_q)$ will have full rank $2n$ at every point of $F_2(M)$.

As we define such a map $f_1$, we must take special care that the map $f_1 \oplus f_1$ extends properly to the closure of $\Sigma_{\tilde g}$ in $M \times M$.  For points $(x,x)$ in the closure of $\Sigma_{\tilde g}$, there exists a sequence of points $(x_m,y_m)$ with corresponding tangent vector $s(x_m,y_m)$.  It will be necessary to study such points separately; this is the content of Section \ref{sec:local}.

Momentarily ignoring this issue near the diagonal, let $\pi : M \times M \to M$ represent the projection onto the first factor.  We make use of the decomposition $s(x,y) = (s_x(x,y),s_y(x,y)) \in T_xM \oplus T_yM$.

We claim that we may replace $g_1$ when the following three conditions are satisfied:
\begin{enumerate}
\setlength\itemsep{-.3em}
\item $\Sigma_{\tilde{g}}$ is a manifold.
\item The restriction of $\pi$ to $\Sigma_{\tilde{g}}$ is an embedding.
\item The section $s_x$ of $TM \big| \pi(\Sigma_{\tilde{g}})$ is nowhere tangent to $\pi(\Sigma_{\tilde{g}}) \subset M$.
\end{enumerate}

\begin{figure}[h!t]
    \centering
        \includegraphics[width=0.4\textwidth]{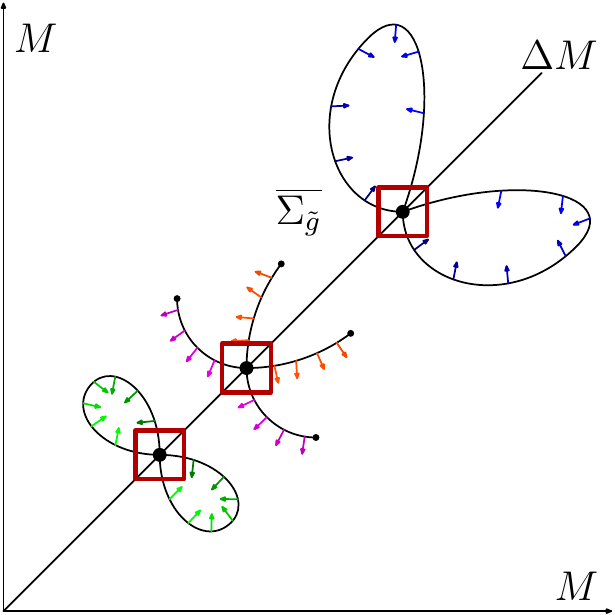} 

        \caption{The $\Z_2$-invariant singularity set $\overline{\Sigma_{\tilde{g}}}$ in $M \times M$, along with its section $s$}
                \label{fig:globalsingproduct}
\end{figure}

\begin{remark}
The key idea of the assumptions is the following: our goal is to define $f_1 \oplus f_1$ carefully on $\Sigma_{\tilde{g}}$.  With the embedding assumption, this singularity set $\Sigma_{\tilde{g}}$ can be associated (via $\pi^{-1}$) with the set $\pi(\Sigma_{\tilde{g}})$ in $M$, and $f_1$ can be explicitly defined in terms of $g_1$ on this set.  (See Figures \ref{fig:globalsingproduct} and \ref{fig:globalsing} in Section \ref{sec:outline} for a depiction.)
\end{remark}

When these conditions are satisfied, we may replace $g_1$ as follows.  Let $f_1$ be $0$ on $\pi(\Sigma_{\tilde{g}})$.  This definition determines $d(f_1)_x$ on vectors tangent to $\pi(\Sigma_{\tilde{g}})$, but this is not an issue, because we only care about the value of $d(f_1)_x$ in the direction of $s_x$, and we have assumed that the section $s_x$ is nowhere tangent.

Now we prescribe $d(f_1)_x$ in the direction of $s_x$.  As $\ker(d\tilde{g})$ is a trivial line bundle on $\pi(\Sigma_{\tilde{g}})$, we may consider a small tubular neighborhood $\pi(\Sigma_{\tilde{g}}) \times [-\delta, \delta]$ extending in the direction of $s_x$, and define on this neighborhood $f_1(x,t) = t(dg_1)_{(x,\pi^{-1}(x))}$; this is well-defined due to injectivity of $\pi$.  Assuming that we may extend $f_1$ to all of $M$, we check that $(f_1 \oplus f_1, \tilde{g})$ is an immersion. The behavior at points $(x,y) \notin \Sigma_{\tilde{g}}$ is irrelevant, and for $(x,y) \in \Sigma_{\tilde{g}}$, we compute
\begin{align*}
d(f_1 \oplus f_1)_{(x,y)}(s) & = d(f_1)_x(s_x) + d(f_1)_y(s_y) = d(g_1)_{(x,y)}(s_x,0) + d(g_1)_{(y,x)}(s_y,0) \\
& = d(g_1)_{(x,y)}(s_x,0) + d(g_1)_{(x,y)}(0,s_y) = d(g_1)_{(x,y)}(s),
\end{align*}
as desired.  Note that here we have used the $\Z_2$-invariance of $\Sigma_{\tilde{g}}$; in particular, $(x,y)$ is an element of $\Sigma_{\tilde{g}}$ with corresponding tangent vector $s(x,y)$ if and only if $(y,x)$ is an element of $\Sigma_{\tilde{g}}$ with tangent vector $s(y,x)$.

We will see that the three assumptions above are satisfied by a generic $\tilde{g}$ when $q > \frac{7}{2} n$; however, we cannot afford to perturb $\tilde{g}$ on $F_2(M)$, since this could disturb the immersion condition on the non-compact set $F_2(M)$.  Thus it is only safe to perturb $\tilde{g}$ away from the diagonal, and we should guarantee the assumptions near the diagonal using a different method.

This leads us to Section \ref{sec:local}, in which we study the local totally nonparallel condition in detail.  Then, in Section \ref{sec:outline} we will refine the outline given above to incorporate the local argument, and also to include a slight technical modification to the assumptions which allows proof in the case $q \geq \frac72$ (necessary to prove Theorem \ref{thm:main} when $n=2$).  In Section \ref{sec:global} we use transversality theorems (proven in Section \ref{sec:transversality}) to show that the assumptions are truly generic in the stated dimensions.  At various stages, we address technical difficulties of patching together the local and the global.

\section{The local problem and the space of singular symmetric bilinear maps}
\label{sec:local}

The weak Whitney theorems state that a generic map $f: M \to \R^{2n+1}$ is an embedding, and that a generic map $f: M \to \R^{2n}$ is an immersion, hence a local embedding.  Then the Whitney trick may be used to resolve double points globally, leading to a proof of the strong Whitney embedding theorem: that there exists an embedding $M \to \R^{2n}$.  We find that the situation for totally nonparallel immersions is analogous.

More precisely, since any immersion is a local embedding, a sufficient condition for the local, zero-order, two-point embedding condition $f(x) \neq f(y)$ is the first-order, one-point immersion condition.  Here, we expect a second-order, one-point condition which is sufficient for the first-order, two-point, local totally nonparallel condition.  Theorem \ref{thm:snbmap} hints at how this condition will manifest.

\subsection{Locally totally nonparallel immersions}

\begin{definition} A map $f : M \to \R^q$ is \emph{locally totally nonparallel at} $x \in M$ if there exists a neighborhood $U \ni x$ such that $f \big|_U$ is totally nonparallel, and $f$ is called \emph{locally totally nonparallel} if every $x$ has such a neighborhood.
\end{definition}

Fix $x \in M$, a neighborhood $U \ni x$, and local coordinates $(x_1,\dots,x_n)$ on $U$.  We consider $df_x$ as an element of  $\Hom(T_xM, \R^q)$ and $\left(\frac{\partial^2 f_i}{\partial x_j \partial x_k}\right)_x$ as an element of $\Hom(T_xM \circ T_xM, \R^q)$, where $\circ$ represents the symmetric product.  When the coordinates are fixed, we will use the notation $\partial^2 f_x$ to represent the element of $\Hom(T_xM \circ T_xM, \R^q)$.

Recall the notion of nonsingularity for symmetric bilinear maps defined in Section \ref{sec:localex}.

\begin{definition}
We say that a smooth immersion $f: M \to \R^q$ is \emph{semifree} if for every point $x \in M$, there exists a neighborhood $U \ni x$ and local coordinates $(x_1,\dots,x_n)$, such that $\partial^2 f_x$ is nonsingular and $\operatorname{Image}(df_x) \cap \operatorname{Image}(\partial^2 f_x) = 0$.
\end{definition}

We first observe that the semifree condition does not depend on local coordinates.  Indeed, a vector in the image of the second derivative, prior to a change in coordinates, may differ from a vector in the image following a change in coordinates, but the difference will be an element in the tangent space, hence the nonsingularity and trivial intersection of images are both preserved under such changes.

An equivalent definition of semifree for an immersion $f : M \to \R^n$ is that the second fundamental form of $f$ is nonsingular, considered as a symmetric bilinear map taking values in the normal bundle.

Examples of semifree maps are the graphs of the maps given in Corollary \ref{cor:quad} above.

\begin{prop}
\label{prop:localtotallynonparallel}
A smooth semifree immersion $f : M \to \R^q$ is locally totally nonparallel.
\end{prop}

\begin{proof}
Consider an immersion $f : M \to \R^q$, $x \in M$, and local coordinates $(x_1,\dots,x_n)$ at $x$.  Assume that $f$ is not locally totally nonparallel at $x$.  Then for every $\varepsilon > 0$, there exist distinct points $y_1, y_2 \in B_\varepsilon(x)$ and nonzero tangent vectors $u_1 \in T_{y_1}M$ and $u_2 \in T_{y_2}M$ such that $df_{y_1}(u_1) = df_{y_2}(u_2)$.  We assume that each $\varepsilon$ is smaller than the injectivity radius of $M$ and small enough so that the tangent bundle over $B_\varepsilon(x)$ is trivial.  Specifically we will consider $u_i \in \R^n$ and $df$ as a map $M \to \Hom(\R^n,\R^q)$.

We may choose each $u_1$ and $u_2$ so that $w \coloneqq df_{y_1}(u_1) = df_{y_2}(u_2)$ is a unit vector in $\R^m$.  In this case, the sequences $\left\{u_1\right\}$ and $\left\{u_2\right\}$ (indexed by some sequence $\varepsilon \to 0$) are bounded away from $0$, and we may restrict to convergent subsequences.  Then there exists $u \in \R^n$ such that $\left\{u_1\right\}$ and $\left\{u_2\right\}$ both converge to $u$, for otherwise the immersion assumption is violated.  Thus we are justified in writing $u_i = u + \delta_i v_i$ for some unit vectors $v_i$.  Similarly, for some unit vectors $h_i$, we are justified in writing $\varepsilon_i h_i$ to represent the tangent vector in $T_xM$ which corresponds via the exponential map to $y_i$.

Now we write the Taylor series of $df : M \to \Hom(\R^n,\R^q)$ at the point $x$:
$$df_{y_i} = df_x + \partial^2f_x(\varepsilon_i h_i) + \mbox{ higher order terms,} \hspace{.25in} i = 1,2.$$
Each term is an element of $\Hom(\R^n,\R^q)$, so that we may act by all terms on the vector $u_i$.  This yields the equations
$$df_{y_i}(u_i) = df_x(u_i) + \partial^2f_x(\varepsilon_i h_i, u_i) + \mbox{ higher order terms,} \hspace{.25in} i = 1,2.$$
Subtracting for $i=2$ from $i=1$ yields
\begin{align*}
0 & = df_x(u_1 - u_2) + \partial^2f_x(\varepsilon_1 h_1, u_1) - \partial^2f_x(\varepsilon_2 h_2, u_2) + \mbox{ h.o.t.} \\
& = df_x(u_1 - u_2) + \partial^2f_x(\varepsilon_1 h_1, u + \delta_1 v_1) - \partial^2f_x(\varepsilon_2 h_2, u + \delta_2 v_2) + \mbox{ h.o.t.} \\
& = df_x(u_1 - u_2) + \partial^2f_x(\varepsilon_1h_1 - \varepsilon_2h_2, u) + \sum_{i=1,2} (-1)^{i+1}  \partial^2f_x(\varepsilon_ih_i,\delta_iv_i) + \mbox{ h.o.t.}
\end{align*}
The terms in the sum are higher order than the first two terms and may be ignored, so we conclude that
$$- df_x(u_1 - u_2) = \partial^2f_x(\varepsilon_1h_1 - \varepsilon_2h_2, u).$$
In case $u_1 = u_2$ for all $\varepsilon$ near $0$, we conclude that $\partial^2f_x$ evaluates to $0$ when applied to the nonzero vectors $\varepsilon_1h_1 - \varepsilon_2h_2$ and $u$, therefore $\partial^2f_x$ is \emph{not} nonsingular as a symmetric bilinear map.  Otherwise we conclude that there is a nonzero vector in the image of both $df_x$ and $\partial^2f_x$.
\end{proof}

The name ``semifree'' is inspired by the following.  A smooth immersion $f : M \to \R^q$ is called \emph{free} if, at every point $x \in M$, the linear space generated by the first and second partial derivatives has the maximum possible number of dimensions $\frac{n^2+3n}{2}$ (see \cite{GromovEliashberg}).  Free maps play an important role in the study of isometric immersions.  As free maps are semifree, they are locally totally nonparallel.

\subsection{The structure of the space of symmetric bilinear maps}

To fully understand semifree immersions, we must understand the nonsingularity condition for symmetric bilinear maps.  Let $n$ and $p$ be natural numbers such that $n < p$.  Recall that a (symmetric) bilinear map $B : \R^n \times \R^n \to \R^p$ is \emph{nonsingular} provided that $B(x,y) = 0$ only when $x = 0$ or $y = 0$.  

Nonsingular symmetric bilinear maps were studied, in part, due to the fact that they induce embeddings of real projective spaces (see \cite{James4}).  Most notably, they were studied in a series of articles by K.Y.\ Lam (e.g.\  \cite{Lam1}, \cite{Lam6}, \cite{Lam4}, \cite{Lam3}, \cite{Lam2}, \cite{Lam5}).  Nonsingular bilinear maps (not necessarily symmetric) are related to immersions of real projective spaces, to Adams' study of vector fields on spheres \cite{Adams}, to Property P \cite{AdamsLaxPhillips}, \cite{AdamsLaxPhillips2}, and to skew flat fibrations \cite{OvsienkoTabachnikov}, \cite{Harrison}.  Nevertheless, the structure of the space of maps which are nonsingular, considered as a subset of the space of all (symmetric) bilinear maps, appears to be unstudied.

I am grateful to Kee Lam for his assistance and comments regarding this subsection, as well as to Daniel Levine for several useful discussions.

We define
\begin{align*}
\Sigma & = \left\{\mbox{symmetric, bilinear maps } B : \R^n \times \R^n \to \R^p \ \big| \ B \mbox{ is not nonsingular} \right\}.
\end{align*}
There is a convenient alternate description of $\Sigma$.  Consider the map $\R^n \times \R^n \to \R^n \circ \R^n$, given by $(x,y) \mapsto x \circ y$.  Let $S \subset \Proj^{\frac{n^2+n}{2}-1}$ be the image of this map projectivized:
\begin{align*}
s : \Proj^{n-1} \times \Proj^{n-1} \to \Proj^{\frac{n^2+n}{2}-1} : (x,y) \mapsto (2x_1y_1,\dots,x_iy_j+x_jy_i,\dots,2x_ny_n), \hspace{.1in} 1 \leq i \leq j \leq n,
\end{align*}
where we have made use of homogeneous coordinates.

A symmetric bilinear map $B : \R^n \times \R^n \to \R^p$ is equivalent to a map $B \in \Hom(\R^n \circ \R^n, \R^p)$.

\begin{lem}
For a symmetric bilinear map $B : \R^n \times \R^n \to \R^p$, $B \in \Sigma$ if and only if $\mathbb{P}(\ker(B)) \cap S \neq \emptyset$.
\label{lem:nonsing}
\end{lem}

\begin{proof} Let $(x,y) \in \R^n \times \R^n$ for nonzero $x,y$.  Then $x \circ y \in S$, and $B(x,y) = B(x \circ y)$, so the forward implication holds.  Conversely, if $t \in S$, then it has a preimage $(x,y) \in R^n \times \R^n$ with $x,y$ nonzero.  Then $x \circ y = t$, and $B(x,y) = B(t)$, so the backward implication holds.
\end{proof}

In this sense we may view $\Sigma$ as a subset of $\Hom(\R^n \circ \R^n, \R^p)$, consisting of maps $B$ for which $\mathbb{P}(\ker(B)) \cap S \neq \emptyset$.

\begin{remark}
Had we forgone the assumption of symmetry, we would instead be considering the Serge variety, i.e.\ the image of the Segre embedding
\begin{align*}
\Proj^{n-1} \times \Proj^{n-1} \to \Proj^{n^2-1} \ \ \ : (x,y) \mapsto (x_1y_1,\dots,x_iy_j,\dots,x_ny_n), \hspace{.1in} 1 \leq i, j \leq n,
\end{align*}
which may equivalently be defined as the projectivization of the set of rank-$1$ $n \times n$ matrices.  The above lemma and the next theorem hold in the non-symmetric situation with essentially the same argument.
\end{remark}

\begin{thm}
\label{thm:strat}
The set $\Sigma$ is a stratified subset of $\Hom(\R^n \circ \R^n, \R^p)$, with codimension $p - 2n + 2$.
\end{thm}

\begin{proof}  Let $\Sigma_k$ be the submanifold of $\Hom(\R^n \circ \R^n, \R^p)$ consisting of maps with kernel dimension $k$.  We will show that $\Sigma_k$ admits a stratification into strata $\Sigma_{kj}$ for which the largest-dimensional stratum is precisely those maps which are nonsingular, assuming any such maps exist.  We can then take all pairs $k,j$ such that $\Sigma_{kj}$ consists of singular maps to obtain a stratification of $\Sigma$.

Let us consider $k$ fixed, let $d = \dim(\R^n \circ \R^n) = \frac{n^2+n}{2}$, and let $\textnormal{Gr}_{k-1}(d-1)$ be the (projective) Grassmann manifold of $(k-1)$-planes in $\Proj^{d-1}$.  Let $E$ be the projectivized tautological bundle over $\textnormal{Gr}_{k-1}(d-1)$, for which the fiber over $P \in \textnormal{Gr}_{k-1}(d-1)$ is $P$ itself.  We consider $E$ as a subset of $\textnormal{Gr}_{k-1}(d-1) \times \Proj^{d-1}$.  Now let $F \subset \textnormal{Gr}_{k-1}(d-1) \times \Proj^{d-1}$ be the bundle obtained by intersecting $E$ with the direct product
\[
\textnormal{Gr}_{k-1}(d-1) \times S \subset \textnormal{Gr}_{k-1}(d-1) \times \Proj^{d-1}.
\]

Then apply Upper Semicontinuity of Fiber Dimension (see e.g.\ \cite{Vakil}, Exercise 18.1.C) to the projection morphism $F \to \textnormal{Gr}_{k-1}(d-1)$.  In particular, we obtain a stratification of $\textnormal{Gr}_{k-1}(d-1)$ which distinguishes elements $P$ by the dimension of the intersection $P \cap S$.  The largest-dimensional stratum consists of those planes $P$ for which $P \cap S$ is empty, if any such planes exist.

Now, there is a bundle $\Sigma_k \to \textnormal{Gr}_{k-1}(d-1) : B \mapsto \ker(B)$, and the fiber over $P \in \textnormal{Gr}_{k-1}(d-1)$ is the set of full rank maps in $\Hom(P^\perp,\R^n)$ (see e.g.\ \cite{GolubitskyGuillemin}, Chapter VI Proposition 1.1).  In particular, this bundle projection is a smooth submersion, and so the stratification of $\textnormal{Gr}_{k-1}(d-1)$ pulls back to a stratification of $\Sigma_k$.  The largest-dimensional stratum consists of those maps $B$ for which $\ker(B) \cap S$ is empty, if any such maps exist.

To show that the largest-dimensional singular stratum in $\Sigma$ has codimension $p-2n+2$, first consider maps $B$ with smallest possible kernel dimension $k = \max\left\{ 0 , d - p \right\}$.  These form an open dense subset of $\Hom(\R^n \circ \R^n, \R^p)$.  The subset of these maps such that $\ker(B) \cap S$ has dimension $0$ is a variety, by the previous argument.  Let $t \in S$ and observe that these maps $B$ for which $t \in \ker(B)$ has codimension $p$.  Taking the union over $t \in S$ gives a variety of codimension $p - 2n + 2$, since $\dim(S) = 2n -2$.
\end{proof}

\subsection{Existence of semifree maps}

\begin{thm} \label{thm:semifreeod} The subset of $C^\infty(M,\R^{4n-1})$ consisting of semifree maps is open and dense.
\end{thm}
\begin{proof}
We will define a subset $\mathscr{S}$ of the $2$-jet bundle $J^2(M,\R^q) \to M$ which is independent of the source $x \in \R^n$ and the target $z \in \R^q$, such that $h \in C^\infty(M,\R^q)$ is semifree if and only if its $2$-jet extension $j^2h$ is disjoint from $\mathscr{S}$.

Let $L = \cup_{i = 1}^n L^i \subset \Hom(\R^n,\R^q)$ consist of those linear maps with corank at least $1$.  The complement $L^0 = L^c$ consists of full rank maps.  Define the singularity subset $\mathscr{L}$ of the trivial bundle $\Hom(\R^n,\R^q) \times \Hom(\R^n \circ \R^n, \R^q) \to \Hom(\R^n,\R^q)$ such that the fiber over an element $P \in L$ is the entire space $\Hom(\R^n \circ \R^n, \R^q)$, and the fiber over $P \in L^0$ consists of those elements $B \in \Hom(\R^n \circ \R^n, \R^q)$ such that the composition of $B$ with the projection $\R^q \to P^\perp$, considered as an element of $\Hom(\R^n \circ \R^n, P^\perp)$, is not nonsingular.

Let $\mathscr{S}$ be the subset of the $2$-jet space for which the fiber over $(x,z) \in \R^n \times \R^q$ is $\mathscr{L}$.  If the $2$-jet extension of a map $h \in C^\infty(M,\R^{4n-1})$ misses $\mathscr{S}$, then $h$ is an immersion, since its $1$-jet lies in $S_0$, and moreover $h$ is semifree, by the condition that its $2$-jet at every point induces a nonsingular symmetric bilinear map whose image intersects the tangent space only trivially.  Since $\mathscr{S}$ is open, so is the space of semifree maps, by the definition of the $C^\infty$ topology.

It remains to compute the codimension of $\mathscr{S}$, which is equal to the codimension of $\mathscr{L}$ as a subset of $\Hom(\R^n,\R^q) \times \Hom(\R^n \circ \R^n, \R^q)$.  The space of symmetric bilinear maps $\R^n \times \R^n \to \R^{q-n}$ which are not nonsingular is stratified with codimension $q - 3n + 2$, by Theorem \ref{thm:strat}, and so $\mathscr{L}$ itself has codimension $q - 3n + 2$.  Therefore, when $q = 4n - 1$, $\mathscr{S}$ has codimension $n+1 > n$, hence by Thom Transversality Theorem is missed by the $2$-jet of a generic $h \in C^\infty(M,\R^{4n-1})$.
\end{proof}

\subsection{The local behavior at non-semifree points}
\label{sec:localbehavior}

We study the typical local behavior near a non-semifree point of a map $\tilde{h} : M \to \R^{q-1}$.  We are most interested in maps $M \to \R^{4n-2}$, when semifreeness is not quite generic but fails only at a finite number of points.  In this case we would like to understand the local behavior of a generic map near such a point.  Our study here is analogous to Whitney's study of maps from an $n$-dimensional manifold into $\R^{2n-1}$ \cite{Whitney}.  Generically, the immersion condition fails at a finite number of points, but with a small perturbation one may ensure a certain quadratic local behavior near such points, so that one can see explicitly the double-point curves of $M$ emanating from these points.  In our case, if a map fails to be semifree, then at least with a small perturbation one may ensure a certain cubic local behavior in the directions where semifreedom fails.

Although we are interested in maps to $\R^{4n-2}$, we make definitions and state results for general target dimension $q-1$, though it is safe to at least assume $q \geq 3n$.

\begin{definition} \label{def:genericsing} We say that a smooth map $\tilde{h} : M \to \R^{q-1}$ has \emph{cubic singularities} if the following conditions hold:
\begin{enumerate}
\item The map $\tilde h$ is an immersion.
\item The map $\tilde h$ fails to be semifree at a finite number of points.
\item At each non-semifree point $x$, there exists a neighborhood $U \ni x$ and local coordinates $(x_1,\dots,x_n)$ such that
\begin{enumerate}
\item at the point $x$,
\[
\frac{\partial^2 \tilde h}{\partial x_1 \partial x_2}= 0,
\]
and the map $(\tilde{h}, x_1x_2) : M \to \R^{4n-1}$ is semifree at $x$ (that is, semifreedom of $\tilde{h}$ fails only in the $x_1x_2$ direction), and
\item at the point $x$, and for both $i = 1$ and $i = 2$, the $(3n-1)$ vectors
\[
\frac{\partial\tilde h}{\partial x_j}, \frac{\partial^2\tilde h}{\partial x_i \partial x_k}, \frac{\partial^3\tilde h}{\partial x_1 \partial x_2 \partial x_\ell} \ \ \ \ \ \mbox{for } 1 \leq j,k,\ell \leq n, \ \ \ k \neq i,
\]
are linearly independent.
\end{enumerate}
\end{enumerate}
\end{definition}

Item 3 simply says that if a map fails to be semifree, then at least it looks cubic in the offending directions.  The advantage of this perspective is that we understand precisely how cubic maps fail to be totally nonparallel: at symmetric points on a curve through the non-semifree point.  For example, $x \mapsto (x,x^3)$ fails to be semifree at the origin, and fails to be totally nonparallel because the tangent vectors at $(x,x^3)$ and $(-x,-x^3)$ are parallel for every $x > 0$.  In this case we say that a \emph{double parallel} occurs at $(x,x^3)$ and $(-x,-x^3)$.  We offer another example.

\begin{figure}[h!t]
    \centering
        \includegraphics[width=0.3\textwidth]{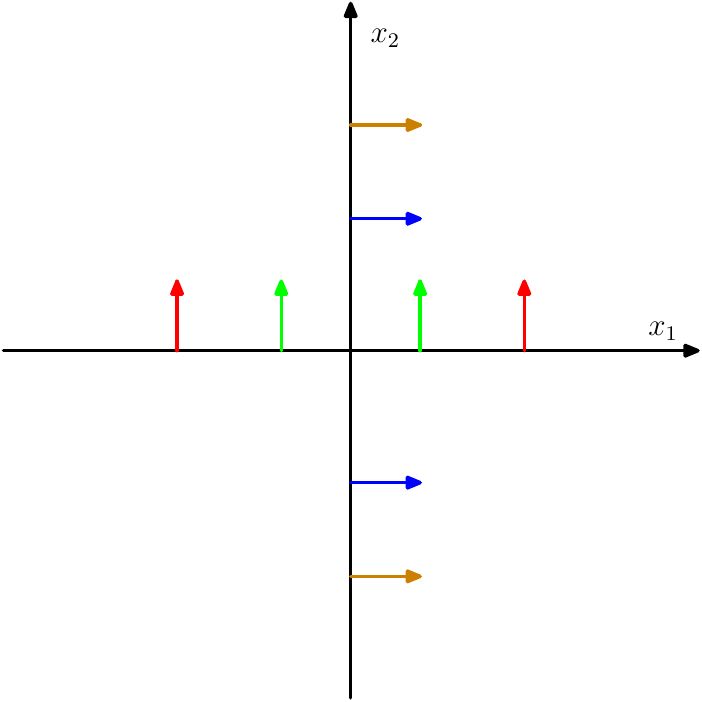} 
        \caption{Local depiction of a non-semifree singularity point in the domain: pairs of identically colored points/tangent vectors get mapped to parallel tangent vectors on the image surface.}
        \label{fig:localsing}
\end{figure}

\begin{figure}[h!t]
    \centering
        \includegraphics[width=0.9\textwidth]{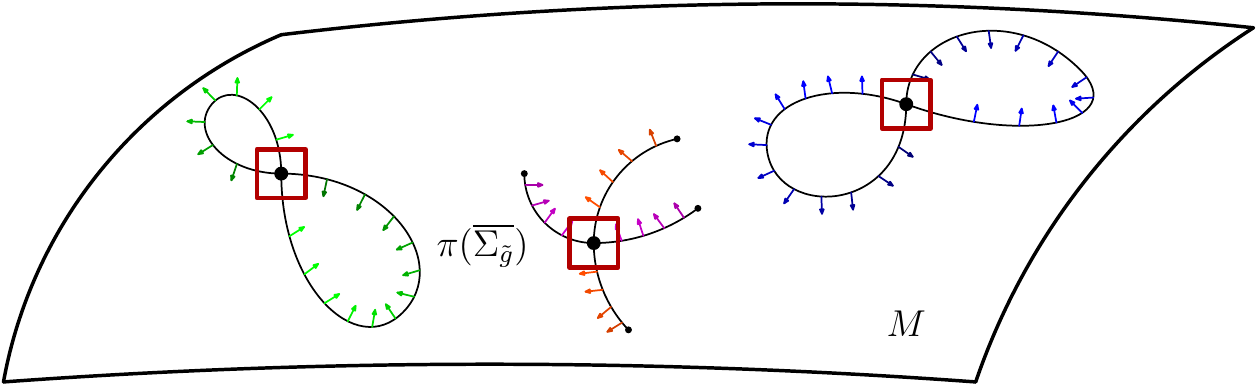} 

        \caption{The projection $\pi(\overline{\Sigma_{\tilde{g}}})$ of Figure \ref{fig:globalsingproduct} to the first factor $M$.  Assuming the sufficient genericity assumption, each dark red square has the local behavior as in Figure \ref{fig:localsing}.  The assumptions in the lemma may be viewed here: $\Sigma_{\tilde{g}}$ is immersed in $M$ along with its nowhere tangent section $s_x$.  In this case the projection $\pi$ is injective, but one could incorporate the more general situation by imagining the upper branch of the middle component of $\pi(\overline{\Sigma_{\tilde{g}}})$ extending further to intersect the right component at some point $x_0$.  In that case, the function $f_i$ would be defined first on the dark red squares around non-semifree points, then in a neighborhood of $x_0$, then extended to all of $\pi(\overline{\Sigma_{\tilde{g}}})$, then extended transversely as determined by the section $s_x$, and finally extended to all of $M$.}
                \label{fig:globalsing}
\end{figure}

\begin{example} \label{ex:cubic} The map $\tilde{h} : \R^3 \to \R^{10}$ given by
\[
\tilde{h}(x_1,x_2,x_3) = (x_1,x_2,x_3,\frac12 x_1^2, \frac12 x_2^2 + x_1x_3, x_2x_3,\frac12 x_3^2,x_1^2x_2,x_1x_2^2,x_1x_2x_3)
\]
is cubically singular.  It fails to be semifree only at the origin $(0,0,0)$, at which semifreedom fails only in the $x_1x_2$ direction (in particular, $\frac{\partial^2 \tilde h}{\partial x_1 \partial x_2} = 0$).  Note that the first seven components of this map were taken directly from the semifree map in Corollary \ref{cor:quad}, except that the component containing $x_1x_2$ is removed; the last three components are the required cubic terms.

It is worth investigating exactly how $\tilde{h}$ fails to be totally nonparallel near the origin.  Let $u_i = \frac{\partial}{\partial x_i}$ be a basis for the tangent plane at an arbitrary point $(x_1,x_2,x_3) \in \R^3$.  Then a vector tangent to the image $\tilde{h}(\R^3)$ at the point $\tilde{h}(x_1,x_2,x_3)$ takes the form

\[
d\tilde{h}_xu = \left( \begin{array}{ccc} 1 & 0 & 0 \\ 0 & 1 & 0 \\ 0 & 0 & 1 \\ x_1 & 0  & 0 \\ x_3 & x_2 & x_1 \\ 0 & x_3 & x_2 \\ 0 & 0 & x_3 \\ 2x_1x_2 & x_1^2 & 0 \\ x_2^2 & 2x_1x_2 & 0 \\ x_2x_3 & x_1x_3 & x_1x_2 \end{array} \right) \left( \begin{array}{c} u_1 \\ u_2 \\ u_3 \end{array} \right) = \left( \begin{array}{cc} u_1 \\ u_2 \\ u_3 \\ x_1u_1 \\ x_3u_1 + x_2u_2 + x_1u_3  \\ x_3u_2 + x_2u_3 \\ x_3u_3 \\ 2x_1x_2u_1 + x_1^2u_2  \\ x_2^2u_1 + 2x_1x_2u_2 \\ x_2x_3u_1 + x_1x_3u_2 + x_1x_2u_3 \end{array} \right).
\]

Assuming that the image of tangent vector $(u_1,u_2,u_3)$ at $(x_1,x_2,x_3)$ is equal to the image of tangent vector $(v_1,v_2,v_3)$ at $(y_1,y_2,y_3)$ yields ten equalities which may be solved.  These can be solved explicitly here (and it may be useful for the reader to do so), but we offer a more general argument instead, which justifies the cubic-singularity definition.  The first three equations relate to the immersion condition and yield $u = v$.  From the fact that $\tilde{h}$ only fails to be semifree in the $x_1x_2$ direction, and from the argument of Proposition \ref{prop:localtotallynonparallel}, we must have either $u = (u_1,0,0)$ or $(0,u_2,0)$.  Without loss of generality we assume $u = (u_1,0,0)$, since the argument is symmetric in $x_1$ and $x_2$.  Then for $u = (u_1,0,0)$, the fact that semifreedom does not fail in other directions (equations 4--7) yields $x_i = y_i$ for $i \neq 2$.  Finally, the linear independence of the cubic partials (equations 8--10) yields $x_i = y_i = 0$ for $i \neq 2$, and $x_2 = -y_2$.

In summary, we find that points $(0,x_2,0)$, $(0,-x_2,0)$ have double parallels (the images of the tangent vector $(u_1,0,0)$) and that points $(x_1,0,0)$, $(-x_1,0,0)$ have double parallels (the images of the tangent vector $(0,u_2,0)$).  Thus the singularity appears in $\R^2$ as two intersecting curves meeting at the non-semifree point (in this case, the $x_1$ and $x_2$ axes which meet at the non-semifree origin); each curve consists of pairs of points with double parallels, and each curve has a vector field in the direction of which the totally nonparallel condition fails (see Figure \ref{fig:localsing}).
\end{example}

\begin{remark} We assume that the non-semifree direction $x_1x_2$ is a ``mixed'' partial instead of a ``pure'' partial because otherwise the vector field, indicating the direction in which the totally nonparallel condition fails, is tangent to the curve along which it fails.  This violates the third assumption in Section \ref{sec:removaltn}.  We observe that this is a reasonable assumption, using an argument analogous to the dimension-count at the end of Theorem \ref{thm:strat}.  Each $t \in \ker(B) \cap S$ corresponds to a nonzero pair $(x,y) \in \R^n \times \R^n$ with $B(x,y) = 0$, and the space of $B$ with $B(x,x) = 0$ for some $x$ has codimension larger than $p - 2n + 2$.
\end{remark}

The example above is instructive because every map with cubic singularities has, in a neighborhood of any non-semifree point, local coordinates such that the singularity appears exactly as in Figure \ref{fig:localsing}.

\begin{definition}
\label{def:sufficientlygeneric} We say that a map $h : M \to \R^{4n-1}$ is \emph{sufficiently generic (with respect to} $i$) if $(h_1,\dots,h_{i-1},h_{i+1},\dots,h_{4n-1})$ has cubic singularities.
\end{definition}


\section{The main argument for Theorem \ref{thm:main}}
\label{sec:outline}

With an understanding of the local totally nonparallel condition, we are now equipped to refine the previous outline of the proof of Theorem \ref{thm:main}.  Recall the goal established in Section \ref{sec:removaltn}, to replace each component $g_i$ of a $\Z_2$-equivariant immersion $g$ with a map of the form $f_i \oplus f_i$.  The following definition should be interpreted in the context of attempting to replace the $i$th component of a map $g$.  Recall Definition \ref{def:sufficientlygeneric} for sufficiently generic maps.

\begin{definition}
We say that a $\Z_2$-equivariant map $g : M \times M \to \R^{4n-1}$ is $i$\emph{-replacement-admissible} if $g$ restricts to an immersion on $F_2(M)$, and if there exists a neighborhood $N$ of the diagonal such that $g \big| _N = (h \oplus h)\big| _N$ for some sufficiently generic (with respect to $i$) semifree map $h : M \to \R^{4n-1}$.
\end{definition}

The name is inspired by the fact that the $i$th component of an $i$-replacement-admissible map is \emph{almost} ready to be replaced by $f_i \oplus f_i$, except that a perturbation may be necessary first.  We will formalize this idea later in the section.  For now, we take a step back to see how this definition fits into our goal of proving Theorem \ref{thm:main}.

\begin{prop}
\label{prop:main}
Let $M$ be a smooth closed $n$-dimensional manifold and let $q = 4n - 1$.  If there exists a $1$-replacement-admissible map $g : M \times M \to \R^{q}$ then there exists a totally nonparallel immersion $M \to \R^q$.
\end{prop}

\begin{remark}
We believe that the above statement holds for $q \geq \frac72 n$.  Most of the argument given here still survives the generalization to such dimensions, except that one needs a more refined understanding of how the cubic-singularity and sufficient-genericity assumptions should be adapted.  This argument, and more generally, an h-principle statement for semifree maps, are current focuses of the author.
\end{remark}

\begin{remark}
We believe that assuming existence of the map $h$ is superfluous (even in the general dimensions of the above remark), but that obtaining a sufficiently generic semifree map $h$ assuming only the existence of a $\Z_2$-equivariant immersion $g : F_2(M) \to \R^q$ is a serious endeavor; see Conjecture \ref{con:main} in Section \ref{sec:intro} and also Section \ref{sec:conclusion} for discussion.
\end{remark}

Theorem \ref{thm:main} follows immediately from combining Proposition \ref{prop:main} with the following:

\begin{prop} Given any smooth closed $n$-dimensional manifold $M$, there exists a $1$-replacement-admissible map $g : M \times M \to \R^{4n-1}$.
\label{prop:admissibleexist}
\end{prop}

\begin{proof}
The space of sufficiently generic semifree maps is open and dense in $C^\infty(M, \R^{4n-1})$.  Any such map $h$ is locally totally nonparallel, so there exists a $\Z_2$-invariant open neighborhood $N$ of $\Delta M$ in $M \times M$ such that $h \oplus h : M \times M \to \R^{4n-1}$ is an immersion on $N - \Delta M$.  In particular (shrinking $N$ if necessary), let $X = M \times M - N$, and we consider $h \oplus h$ as a smooth map on $X / \Z_2$ which is an immersion in a neighborhood $U$ of $\partial N / \Z_2$.  We would like to extend $h \oplus h$ to an immersion $g : X / \Z_2 \to \R^{4n-1}$ such that $g$ is equal to $h \oplus h$ on $U$.

By Whitney \cite{Whitney}, Theorem 7, such an extension is possible provided that $h \oplus h$ is an embedding on $U$ and that the number of singular points (i.e.\ points of nonimmersion) of $h \oplus h$ on $X/\Z_2$ is even.  Let $\Sigma$ represent the singular points of $h \oplus h$; these correspond to (unordered) pairs of points at which the totally nonparallel condition fails for $h$.  By the proof of Theorem \ref{thm:gen} and genericity of $h$, the set of singular points of $h \oplus h$ is zero-dimensional.  Hence there is a well-defined element $[\Sigma] \in H^{2n}(X/\Z_2; \Z_2)$ dual to $\Sigma$.  By the Thom-Porteous formula (\cite{ArnoldEtAl} Chapter 4, 1.4), $[\Sigma] = \bar{w}_{2n}(X/\Z_2)$, which is zero since every $2n$-dimensional manifold immerses in $\R^{4n-1}$.  Hence the number of singular points is indeed even.

It remains to show that there exists an $h$ and $N$ such that $h \oplus h$ is an embedding on $\partial N / \Z_2$, since then it is an embedding on some neighborhood $U$, as desired.  As $\partial N / \Z_2$ is $(2n-1)$-dimensional, the condition of injectivity on $\partial N / \Z_2$ is generic \emph{among maps in} $C^\infty(X/\Z_2, \R^{4n-1})$.  We use Corollary \ref{cor:ttb3} to show that the same holds for maps of the form $h \oplus h$.

In the statement of the corollary, let $X = (M \times M - N')$, where $N'$ is an open neighborhood of the diagonal $\Delta M$ smaller than $N$, and let the singularity set $V \subset J_2^0(X/\Z_2, \R^{4n-1})$ be chosen to represent the ``non-injectivity condition'' on $\partial N / \Z_2$. In particular, $V$ may be defined fiberwise, so that the fiber over $([(x,y)],[(x',y')]) \in X / \Z_2 \times X / \Z_2 - \Delta(X / \Z_2)$ is equal to:
\begin{align*}
\Delta \R^{4n-1} \subset \R^{4n-1} \times \R^{4n-1}  \hspace{.5in} & \mbox{for } ([(x,y)],[(x',y')]) \in \partial N / \Z_2 \times \partial N / \Z_2 - \Delta(\partial N / \Z_2) \\
\emptyset \hspace{1.35in} & \mbox{otherwise.}
\end{align*}
With this definition, $(j_2^0(h \oplus h))$ intersects $V$ if and only if $(h \oplus h)|_{\partial N / \Z_2}$ is not injective.

Now $V$ has (fiber) codimension $4n-1$, and so by Corollary \ref{cor:ttb3}, there is a residual set of $h$ such that $V$ is missed by the $(4n-2)$-dimensional image of the $2$-fold $0$-multijet of $(h \oplus h)|_{\partial N / \Z_2}$.  Since the set of sufficiently generic semifree maps is open and dense, and the set of those satisfying the non-injectivity condition is dense, there exists an element $h$ in the intersection, completing the proof.
\end{proof}

The remainder of the paper is devoted to the proof of Proposition \ref{prop:main}.  An essential portion of the argument is the following lemma, which refines the proof outline of Section \ref{sec:removaltn}.  The idea is to inductively replace each component $g_i$ of some $i$-replacement-admissible map $g = (g_1,\dots,g_q)$ with $f_i \oplus f_i$ while preserving the immersion condition on $F_2(M)$.

Assume that we are on step $i$ of this proposed replacement; in particular, that we have an $i$-replacement-admissible map $g$ such that $\tilde{g} = (F \oplus F, G)$ for some $F \coloneqq (f_1,\dots,f_{i-1}) : M \to \R^{i-1}$ and $G \coloneqq (g_{i+1},\dots g_q): M \times M \to_{\Z_2} R^{q-i}$.  As in Section \ref{sec:removaltn}, define
\[
\Sigma_{\tilde g} \coloneqq \left\{ (x,y) \in F_2(M) \ \big| \ \rank(d\tilde g_{(x,y)}) = 2n - 1 \right\}.
\]

\begin{lem}
\label{lem:assumption}  Let $M$ be a smooth closed $n$-dimensional manifold, and let $g : M \times M \to \R^q$ be $i$-replacement-admissible.  Let $\pi : M \times M \to M$ represent the projection onto the first factor.  Suppose that
\begin{enumerate}
\setlength\itemsep{-.3em}
\item $\Sigma_{\tilde{g}}$ is a manifold.
\item The restriction of $\pi$ to $\Sigma_{\tilde{g}}$ is an immersion.
\item The section $s_x$ of $TM \big| \pi(\Sigma_{\tilde{g}})$ is nowhere tangent to $\pi(\Sigma_{\tilde{g}}) \subset M$.
\item The restriction of $\pi$ to $\Sigma_{\tilde{g}}$ has no triple points.
\item There exists $\beta > 0$ such that the restriction of $\pi$ to the closure $\overline{\Sigma_{\tilde{g}}}$ is injective, except that there may exist (finitely many) pairs of points $(x,y)$, $(x,y')$ in $\Sigma_{\tilde{g}}$, at least distance $\beta$ apart, and both at least distance $\beta$ from the diagonal $\Delta M$.  At each such pair, the projected kernel vectors $s_x(x,y)$, $s_x(x,y')$ in $T_xM$ are not parallel.
\end{enumerate}
Then there exists a function $f_i : M \to \R$ such that $(f_i \oplus f_i,\tilde{g})$ is an immersion on $F_2(M)$ and $(f_i,\tilde{h})$ is semifree.

Moreover, in the case $q = 4n - 1$, there exists a perturbation of $(f_1,\dots,f_i)$ in $C^\infty(M, \R^i)$ such that the resulting function $(f_1,\dots,f_i,h_{i+2},\dots,h_{4n-1})$ has cubic singularities.  In particular, $(f_1 \oplus f_1, \dots, f_i \oplus f_i, g_{i+1}, \dots, g_{4n-1})$ is $(i+1)$-replacement-admissible.
\end{lem}

\begin{proof}
We begin at the finitely many points where $\tilde{h}$ fails to be semifree.  Since $\tilde{h}$ is assumed to have cubic singularities, we understand the local behavior of $\tilde{h}$ at these non-semifree points.  At any such point we may choose coordinates $(t_1,t_2)$ as in Definition \ref{def:genericsing}, so that $\tilde{h}$ fails to be semifree only with respect to the mixed partial direction $t_1t_2$, and since $h$ itself is semifree, it must be the case that $\frac{\partial h_i}{\partial t_1 \partial t_2} \neq 0$ at the non-semifree point.  This is the condition we must preserve when defining $f_i$, since we require that $(f_i, \tilde{h})$ is still semifree.

Consider a small two-dimensional square $[-\delta, \delta] \times [-\delta, \delta]$ in these coordinates, so that the non-semifree point occurs at $(0,0)$.  On this square define $f_i(t_1,t_2) = \frac{t_1t_2}{\delta^2}$.  This guarantees that the new map $(f_i, \tilde{h})$ is semifree at these points.

Now we would like to extend $f_i$ to the rest of $\Sigma_{\tilde{g}}$.  We assume momentarily that the restriction of $\pi$ to $\overline{\Sigma_{\tilde{g}}}$ is injective.  In this case the proof outline of Section \ref{sec:removaltn} goes through almost unchanged.  In particular, we let $f_i$ be $0$ on $\pi(\overline{\Sigma_{\tilde{g}}})$, and note that this is consistent with the definition of $f_i$ on the $2$-dimensional neighborhoods around non-semifree points.  We may again prescribe $d(f_i)_x$ in the direction of $s_x$, using exactly the argument in the outline of Section \ref{sec:removaltn}.  In this case, we have defined $f_i$ on a closed set, and we may extend $f_i$ to all of $M$, and the resulting map satisfies the conditions that $(f_i \oplus f_i, \tilde{g})$ is an immersion and that $(f_i, \tilde{h})$ is semifree.

We adapt the proof in case we do not necessarily have injectivity of $\pi$ on $\overline{\Sigma_{\tilde{g}}}$.  In this case, instead of defining $f_i$ as identically zero on $\Sigma_{\tilde{g}}$, we may define $f_i$ as zero on the finitely many self-intersection points of $\pi(\Sigma_{\tilde{g}})$ and extend $f_i$ in a $2$-dimensional square determined by the \emph{two} transverse pieces of the section $s_x$.  Note that this may determine $f_i$ along some arc of $\pi(\Sigma_{\tilde{g}})$; for example, when $M$ is $2$-dimensional, this determines $f_i$ in a neighborhood of $x$.  However, we reiterate that the value of $f_i$ is irrelevant; the only important value is that of the derivative in the direction of the section $s_x$, which we may still prescribe by the transversality assumptions.

Next, we argue in the case $q = 4n - 1$ that we may perturb $(F,f_i)$ to guarantee cubic singularities of $\tilde{h}_{i+1} \coloneqq (F,f_{i},h_{i+2},\dots,h_{4n-1})$.  Note that such singularities are guaranteed with an arbitrarily small perturbation of $\tilde{h}_{i+1}$; however, such a perturbation is dangerous, because it may disrupt $\tilde{h}_{i+1} \oplus \tilde{h}_{i+1}$ on the boundary of $N$, and it is not trivial to re-glue this perturbed map to a new immersion $g$ (though we suspect that one could make an argument using this method).  Thus we would like to only perturb the functions $(F,f_i)$ which we have already replaced, since then there is no gluing problem to worry about.

We are interested in the set of non-semifree points of $\tilde{h}_{i+1}$; we would like to ensure that this is a finite set with cubic singularities.  Consider its superset
\[
\Sigma = \left\{ x \in M \ \big| \ (h_{i+2},\dots,h_{4n-1}) : M \to \R^{4n-i-2} \mbox{ is not semifree} \right\}.
\]
Let $H \coloneqq (h_{i+2},\dots,h_{4n-1})$.  By genericity of $H$, and by Theorem \ref{thm:strat}, $\Sigma$ is stratified with $\codim(\Sigma) = \max\left\{0,n-i\right\}$, or dimension $\min \left\{n, i \right\}$.  We will define a singularity set $\mathscr{S}$ in the $2$-jet space $J^2(M, \R^i)$ (thought of as the $2$-jet space corresponding to the functions $(F,f_i)$) over the points $\Sigma$ as follows.

Let us assume that $i < 2n -1$ (the other case is similar).  Then by genericity $H$ is an immersion, and its $2$-jet induces a section $B_H$ of $\Hom(TM \circ TM, \R^{4n - i - 2} / TM)$.   By Lemma \ref{lem:nonsing}, $H$ fails to be semifree at $x$ if and only if $\ker(B_H(x))$ intersects $S \subset T_xM \circ T_xM$.  Let $S_x^\prime$ be the intersection at such a point.  Then the singularity set $\mathscr{S}$ in the $2$-jet space $J^2(M, \R^i)$ may be defined fiberwise over $x \in \Sigma$ as the set of maps $B \in \Hom(T_xM \circ T_x M, \R^i)$ such that $\ker(B) \cap S_x^\prime$ is nonempty, and with empty fiber over $x \notin \Sigma$.  With this definition, the map $(F,f_i,H)$ is semifree if and only if the $2$-jet of $(F,f_i)$ misses $\mathscr{S}$.  This is analogous to Lemma \ref{lem:nonsing} with $S_x^\prime$ in place of $S$.

It remains to compute the codimension of $\mathscr{S}$.  For fixed $t \in S_x^\prime$, the codimension of the set of $B$ with $t \in \ker(B)$ is $i$.   The expected dimension of $S_x^\prime$ is $\max\left\{0,i-n\right\}$, since $S$ has dimension $2n-2$ and almost all $B$ have a kernel of codimension $3n-i-2$.  Thus taking the union over the $(\max\left\{0,i-n\right\})$-dimensional $S_x^\prime$ gives a singularity set $\mathscr{S}$ of codimension $\min\left\{i,n\right\}$.
Then $\mathscr{S}$ is intersected by the $2$-jet of a generic $(F,f_i)$ on a set of codimension $\min\left\{i,n\right\}$, which has dimension $0$ in $\Sigma$.  Thus for generic $(F,f_i)$ we have a finite number of non-semifree points of the new $\tilde{h}_{i+1}$.

Finally, we show that we may guarantee cubic singularities with a small perturbation of $\tilde{h}_{i+1}$, being careful not to disrupt the gluing on the boundary of $N$.  Around each of these finitely many points $x_j$ we may consider a small neighborhood $U_j$, so that for some neighborhood $N$ of the diagonal, no point $(x,y) \in U_j \times U_k$ lies in the boundary of $N$ (so that a perturbation of $\tilde{h}_{i+1}$ on the $U_j$ does not disrupt the gluing to $g$ on the boundary of $N$).  Then by a small perturbation of $\tilde{h}_{i+1}$ on the $U_j$, we may guarantee cubic singularities.
\end{proof}

\begin{example}  Suppose we start with $h : \R^2 \to \R^7$, $h(x,y) = (xy,x^2-y^2,xy,x^2-y^2,xy,x,y)$.  Since the last four component functions guarantee semifreedom on their own, the first three replacements may very well yield $(0,0,0,x^2-y^2,xy,x,y)$.  After the third replacement, we examine the set $\Sigma$ at which $(h_5,h_6,h_7)$ fails to be semifree, which is the entire $\R^2$.  In this case, $B_H(x) : T_xM \circ T_xM \to \R^3 / T_xM \simeq \R$ maps $u \circ v$ to $u_1v_2+u_2v_1$.  The kernel is $2$-dimensional, spanned by $e_1 \circ e_1$ and $e_2 \circ e_2$, or $1$-dimensional in the $2$-dimensional projectivization of $T_xM \circ T_xM$.  On the other hand, $S$ is $2$-dimensional in the same projectivization, and their intersection $S_x^\prime$ is $1$-dimensional.  The singularity set $\mathscr{S}$ is defined fiberwise in $J^2(M,\R^3)$ over $x \in \Sigma = \R^2$ as $S_x^\prime$.  By the argument above, the set of $B \in \Hom(\R^2 \circ \R^2, \R^3)$ such that $\ker(B) \cap S_x^\prime \neq \emptyset$ has codimension $2$.  Therefore for generic $(f_1,f_2,f_3)$, the preimage $(j^2(f_1,f_2,f_3))^{-1}(\mathscr{S})$ has codimension $2$, or dimension $0$ in $\R^2$.
\end{example}

\subsection{The inductive application of Lemma \ref{lem:assumption} to prove Proposition \ref{prop:main}}

Proposition \ref{prop:main} may now be argued inductively.  Assume at step $i$ that there exists an $i$-replacement-admissible map $g = (F \oplus F,g_i, G)$ for some $F \coloneqq (f_1,\dots,f_{i-1}) : M \to \R^{i-1}$ and $G \coloneqq (g_{i+1},\dots g_{4n-1}): M \times M \to_{\Z_2} R^{4n-i-1}$.  Existence in the base case $i = 1$ is guaranteed by Proposition \ref{prop:admissibleexist}.  Proposition \ref{prop:main} will follow once we argue that we may guarantee the assumptions of Lemma \ref{lem:assumption} with an arbitrary small perturbation of $\tilde{g}$ for any $i$-replacement-admissible map $g$, since the upshot of Lemma \ref{lem:assumption} leaves us with an $(i+1)$-replacement-admissible map of the desired form $(F \oplus F, f_i \oplus f_i, G)$.  To formalize this perturbation of $\tilde{g}$, we define the relevant function space for fixed $i$, $N$, and $\tilde{h} : M \to \R^{q-1}$:
\begin{align*}
\mathscr{H}_N^i(\tilde{h})  \coloneqq \Big\{ \tilde{g} = (F \oplus F, G) : M \times M \to_{\Z_2} \R^{q-1}
\ \Big| \ \tilde{g} |_N = \tilde{h} \oplus \tilde{h}|_N \Big\}.
\end{align*}
Since the set $\Sigma_{\tilde{g}}$ depends only on $\tilde g$, and not $g_i$, we study the assumptions of Lemma \ref{lem:assumption} by studying the spaces $\mathscr{H}_N^i(\tilde{h})$.

We digress with a quick dimension-count to ensure that a few of the assumptions of Lemma \ref{lem:assumption} are expected to hold when $q \geq \frac72 n$.  Informally, $\Sigma_{\tilde{g}}$ is the space where $\tilde{g} : M \times M \to \R^{q-1}$ has corank one, so it is generically a manifold with expected codimension $q-2n$ and expected dimension in $M \times M$ equal to $4n-q$.  Thus a generic map from $\Sigma_{\tilde{g}}$ to the $n$-dimensional manifold $M$ is an immersion when $2(4n-q) \leq n$, or $\frac72 n \leq q$, and has no triple points as long as $\frac{10}{3} n < q$, and is an embedding as long as $\frac72 n < q$.   The other assumptions are more technical and will be discussed rigorously in future sections.

\begin{proof}[Proof of Proposition \ref{prop:main}] \ \ Though structured as an outline, the following proof is rigorous besides its dependence on results of future sections.  The assumptions of Lemma \ref{lem:assumption} are guaranteed in four steps, roughly summarized as follows: that assumptions 1--5 hold on $N$ by a small perturbation of $h$, that assumptions 1--3 hold outside $N$ from a perturbation of $\tilde{g}$ on $\mathscr{H}_N^i(\tilde{h})$, that assumptions 4--5 hold outside $N$ for points within distance $\beta$, and that assumptions 4--5 hold outside $N$ for points separated by distance $\beta$. \\

\noindent \emph{Step 1.}  Given a sufficiently generic semifree map $h$, there exists an arbitrary small perturbation (we will also call $h$) and open neighborhoods $N \supset N_\alpha \supset \Delta(M)$ such that:
\begin{enumerate}
\setcounter{enumi}{-1}
\setlength\itemsep{-.3em}
\item the restriction $(h \oplus h) \big|_N$ is a $\Z_2$-equivariant immersion,
\item the set $\Sigma_{\tilde{g}} \cap N$ is a manifold,
\item the restriction of $\pi$ to $\Sigma_{\tilde{g}} \cap N$ is an immersion,
\item the section $s_x$ is nowhere tangent to $\pi(\Sigma_{\tilde{g}} \cap N)$,
\item there are no pairs of points $(x,y), (x,y') \in \overline{\Sigma_{\tilde{g}}} \cap N$, and
\item there are no pairs $(x,y) \in \Sigma_{\tilde{g}} \cap N^c$, $(x,y') \in \overline{\Sigma_{\tilde{g}}} \cap N_\alpha$.
\end{enumerate}

The zeroth item is guaranteed by semifreedom.  The other five loosely correspond to the five assumptions of Lemma \ref{lem:assumption} and may be guaranteed as follows, assuming a small perturbation near non-semifree points of $\tilde{h}$.  In particular, there are finitely many non-semifree points of $\tilde{h}$ by sufficient genericity, and as in the end of the proof of Lemma \ref{lem:assumption}, we may choose small neighborhoods around each such point and perturb $\tilde{h}$ only on these neighborhoods (so that the gluing of $h \oplus h$ to $g$ is not disrupted).  Then, the first item is guaranteed by the cubic singularities condition, since we have an explicit understanding of $\Sigma_{\tilde{g}}$ in a neighborhood of a non-semifree point.  Additionally, we may arrange so that the directions $x_1$ and $x_2$ of non-semifreedom are not parallel to the first factor $M$, assuring the second and fourth items.  Item 3 follows from the definition of cubic singularities, since the direction of failed semifreedom is mixed, not pure (see the remark following Example \ref{ex:cubic}).  Finally, $\pi(\Sigma_{\tilde{g}} \cap N^c)$ is fixed, and may be avoided by perturbing the finitely-many points of non-semifreedom. \\

\noindent \emph{Step 2.}  Next we use certain transversality theorems, formulated specifically for the function space $\mathscr{H}_N^i(\tilde{h})$, to show that (for $q \geq \frac72 n$) there exists an open dense subset of $\mathscr{H}_N^i(\tilde{h})$ consisting of functions such that $\Sigma_{\tilde{g}} \cap N^c$ is a manifold (Lemma \ref{lem:manifold}), the restriction of $\pi$ to $\Sigma_{\tilde{g}} \cap N^c$ is an immersion, and the section $s_x$ is nowhere tangent to $\pi(\Sigma_{\tilde{g}} \cap N^c)$ (Lemma \ref{lem:immersion}).  Note that a perturbation of $\tilde{g}$ does not disrupt the properties ensured in Step 1, except that we may need to use a slightly smaller neighborhood $N_\alpha$. \\

\noindent \emph{Step 3.}  A classical fact about immersions states that not only are immersions injective on small neighborhoods, but that there exists a number $r$ such that all arbitrarily small perturbations of an immersion are injective on neighborhoods of size $r$ (see \cite{GolubitskyGuillemin}, II.5.Lemma A).  Using a similar argument, we show in Lemma \ref{lem:taylor} that since $\pi$ is an immersion of $\Sigma_{\tilde{g}}$, there exists a small number $\beta$ such that $\pi$ is injective on sets of size $\beta$, and also that this injectivity holds for arbitrarily small perturbations of $\tilde{g}$. This step does not itself require a perturbation, but it sets us up for the final one. \\

\noindent \emph{Step 4.}  In Section \ref{sec:ttb}, we show that we may make one final perturbation to obtain that (for $q \geq \frac72 n$) there are finitely many $(x,y), (x,y') \in \overline{\Sigma_{\tilde{g}}} \cap N_\alpha^c$, and that any such pairs satisfy $\dist((x,y),(x,y')) > \beta$ and that the vectors $s_x(x,y)$ and $s_x(x,y')$ are not parallel.  In particular, we assure this property with a perturbation within the space $\mathscr{H}_{N_{\min\left\{\alpha,\beta\right\}}}^i(\tilde{h})$, where $\alpha$ and $\beta$ have been defined in Steps 1 and 3, respectively.  In this case, we might change $\tilde{g}$ on $N - N_{\min\left\{\alpha,\beta\right\}}$, but all of the previously established conditions are open for $\tilde{g}$ on $N_\alpha^c$ so are still satisfied. \\

Hence, assuming the results of future sections, the proof of Proposition \ref{prop:main} is complete.
\end{proof}

\section{The technical assumptions of Lemma \ref{lem:assumption}}
\label{sec:global}

\subsection{Transversality theorems}
\label{sec:transstate}

In this section we carry out Steps 2--4 of the proof of Proposition \ref{prop:main}; in particular, we verify that the technical assumptions of Lemma \ref{lem:assumption} can be assured by small perturbations provided that $q \geq \frac72 n$.

At stage $i$ of our replacement, we have an $i$-replacement-admissible map $g = (F \oplus F, g_i, G)$, and we aim to replace $g_i$ with a suitable map $f_i \oplus f_i$.  For the remainder of this paper, we may consider $i \in \left\{ 1, \dots, q \right\}$, $N$, and $\tilde{h} : M \to \R^{q-1}$ fixed, and we study the space $\mathscr{H}_N^i(\tilde{h})$ defined in the previous section:
\begin{align*}
\mathscr{H}_N^i(\tilde{h})  \coloneqq \Big\{ \tilde{g} = (F \oplus F, G) : M \times M \to_{\Z_2} \R^{q-1}
\ \Big| \ \tilde{g} |_N = \tilde{h} \oplus \tilde{h}|_N \Big\}.
\end{align*}
Although $\tilde{g}$ is an element of $C^\infty(M \times M, \R^{q-1})$, it takes a certain special form, due to the restriction on the first $i-1$ coordinate functions and the $\Z_2$-equivariance.  Therefore, the jet extensions of such $\tilde{g}$ do not occupy the entire jet spaces $J^s(M\times M, \R^{q-1})$.  (For example, for a function $f : \R \to \R$, the second derivative of $f \oplus f$ has no mixed partials).  Let $\mathscr{J}$ be the subbundle of $J^s(M\times M, \R^{q-1})$ generated by $s$-jets of functions of the form $(F \oplus F, G)$.  In particular, $\mathscr{J}$ is the image of the natural inclusion map $J^s(M, \R^{i-1}) \times J^s(M, \R^{i-1}) \times J^s(M \times M, \R^{q-i}) \to J^s(M\times M, \R^{q-1})$.

\begin{thmn}[Transversality Theorem A] Let $M$ be a compact manifold, $V \subset J^s(M \times M, \R^{q-1})$ be a stratified set such that $V \pitchfork \mathscr{J}$, and let $N$ be a fixed open neighborhood of $\Delta M$.  Then
\begin{align*}
\left\{ \tilde{g} = (F \oplus F, G) \in \mathscr{H}_N^i(\tilde{h}) \ \Big| \ j^s(\tilde{g}) |_{(M \times M) \setminus N} \pitchfork V \right\}.
\end{align*}
is open and dense in $\mathscr{H}_N^i(\tilde{h})$.
\end{thmn}

Briefly, the result states that we may apply ordinary transversality theorems for functions $\tilde{g}$ to the subset of those functions which are of the form $(F \oplus F, G)$, \emph{provided that the singularity set $V$ is transverse to the space of $s$-jets of such functions}.  The necessity of this assumption may be seen by simply taking $V = \mathscr{J}$; in this case the upshot is clearly false.  However, this important assumption was not mentioned in \cite{szucs}, where similar transversality statements were claimed and used.  It is also worth mentioning that although our function $\tilde{g} \in \mathscr{H}_N^i(\tilde{h})$ takes a special form on $N$, the statement claims nothing about transversality on $N$.

We will use Transversality Theorem A to complete Step 2 of the outline of the proof of Theorem \ref{thm:main}.  To complete Step 4 we require a second transversality theorem, which may be considered as a multijet version of Transversality Theorem A.  Let $N' \subset N$ be a $\Z_2$-invariant open neighborhood of the diagonal.  Let $X = M \times M - N'$.  We will study the $2$-fold $s$-multijet of $\tilde{g}|_X$:
\begin{align*}
j_2^s\tilde{g} & : X \times X - \Delta X \to \R^{q-1} \times \R^{q-1} : (x,y,x',y') \mapsto (j^s\tilde{g}(x,y),j^s\tilde{g}(x',y')).
\end{align*}

Now $j_2^s\tilde{g}$ is a section of the multijet bundle $J_2^s(X, \R^{q-1})$.  Let $\mathscr{J}_2$ be the subbundle generated by such jets.

\begin{thmn}[Transversality Theorem B] Let $M$ be a compact manifold, $N' \subset N$ fixed $\Z_2$-invariant open neighborhoods of the diagonal $\Delta M$, $X = M \times M - N'$, $Z$ a fixed open neighborhood of the diagonal $\Delta X$, and $V \subset J_2^s(X,\R^{q-1})$ a stratified set such that $V \pitchfork\mathscr{J}_2$.  Then
\begin{align*}
\left\{ \tilde{g} \in \mathscr{H}_N^i(\tilde{h}) \ \Big| \ j_2^s(\tilde{g}) |_{Z^c} \pitchfork V \right\}.
\end{align*}
is open and dense in $\mathscr{H}_N^i(\tilde{h})$.
\end{thmn}

We postpone the proofs of the Transversality Theorems until Section \ref{sec:transversality}, opting first to formalize the proof of Proposition \ref{prop:main} assuming the transversality statements.

\subsection{Applications of Transversality Theorem A to Step 2}
\label{sec:step2}

Here we apply Transversality Theorem A to formalize Step 2 of the proof of Proposition \ref{prop:main}.  The arguments of this section are similar in flavor to the singularity-theoretic arguments of Sections \ref{sec:review}--\ref{sec:removaltn}; we use the notation introduced there.

Let $N^c$ represent the complement of $N$ in $M \times M$.  Referring to the first three assumptions in Lemma \ref{lem:assumption} for guidance, we define the following subsets of $\mathscr{H}_N^i(\tilde{h})$:
\begin{align*}
\hman & = \big\{ \tilde{g} \in \mathscr{H}_N^i(\tilde{h}) \ \big| \ \Sigma_{\tilde{g}} \cap N^c \mbox{ is a manifold} \big\} \\
\himm & = \big\{ \tilde{g} \in \hman \ \big| \ \pi |_{\Sigma_{\tilde{g}} \cap N^c} \mbox{ is an immersion} \big\} \\
\hnot & = \big\{ \tilde{g} \in \hman \ \big| \ s_x  \mbox{ is nowhere tangent to } \pi(\Sigma_{\tilde{g}} \cap N^c) \subset M \big\}. 
\end{align*}
\begin{lem} 
\label{lem:manifold}
The subset $\mathscr{H}_{\mbox{\emph{man}}}$ is open and dense in $\mathscr{H}_N^i(\tilde{h})$. 
\end{lem}

\begin{proof}
We associate each fiber of the $1$-jet bundle $J^1(M \times M, \R^{q-1})$ with $\Hom(\R^n \times \R^n, \R^{q-1})$.  Let $S_r$ be the subset which fiberwise consists of linear maps with corank $r$, and let $V$ be the union of $S_r$ for positive $r$.  Then we may write
\[
\Sigma_{\tilde{g}} \cap N^c = \left\{ (x,y) \in N^c \ \big| \ j^1\tilde{g}(x,y) \in V \right\}.
\]
Then $V$ is a stratified (by corank) subset of $J^1(M \times M, \R^{q-1})$, and moreover, $V \pitchfork \mathscr{J}$, since in this case $\mathscr{J} = J^1(M \times M, \R^{q-1})$.  Hence Transversality Theorem A applies, and so the set
\[
\left\{ \tilde{g} \in \mathscr{H}_N^i \ \big| \ j^1\tilde{g}|_{N^c} \pitchfork V \right\}
\]
is open and dense in $\mathscr{H}_N^i$.  For such $\tilde{g}$, the set $\Sigma_{\tilde{g}} \cap N^c = (j^1\tilde{g}|_{N^c})^{-1}(V)$ is a manifold by transversality. 
\end{proof}

For the next result, we require the notion of \emph{intrinsic derivative} due to Porteous. We rely heavily on the treatment of Golubitsky and Guillemin \cite{GolubitskyGuillemin}.

Let $E \to X$ and $F \to X$ be vector bundles, and let $\rho : E \to F$ be a bundle map, considered as a map $X \to \Hom(E,F)$.  Assume that for some point $x$, $\rho(x)$ has corank $r$.  Then the intrinsic derivative is defined as the map
$$ T_xX \xrightarrow{(d\rho)_x} T_{\rho(x)} \Hom(E,F) \simeq \Hom(E,F) \longrightarrow{} \Hom(\Ker(\rho(x)), \Coker(\rho(x))),$$
where the last arrow is given by restricting and projecting.  It is shown in \cite{GolubitskyGuillemin} that this map does not depend on the choice of trivializations.

Let us examine the special case when $\rho$ is the derivative of some map.  Consider $\phi : X \to Y$, let $E = TX$ and $F = \phi^*(TY)$ be the pullback bundle.  Let $\rho : E \to F$ be the map $d\phi$.  Then the intrinsic derivative is the map
\[
D(d\phi)_x : T_xX \to \Hom(\Ker((d\phi)_x), \Coker((d\phi)_x)),
\]
which (\cite{GolubitskyGuillemin}, p. 151, Ex 1) can be viewed as a map
\[
D(d\phi) \in \Hom(TX,\Hom(\Ker(d\phi),\Coker(d\phi))).
\]

Returning to our setup, in particular the previous proof, we associate each fiber of the $1$-jet bundle $J^1(M \times M, \R^{q-1})$ with $\Hom(\R^n \times \R^n, \R^{q-1})$.  Let $L^r \subset \Hom(\R^n \times \R^n, \R^q)$ consist of those linear maps with corank $r$.  Let $S_r \subset J^1(M \times M, \R^{q-1})$ be the submanifold which is fiberwise $L_r$.  The argument of Lemma \ref{lem:manifold} applies to each individual $S_r$, and the corank formula gives $\codim (S_r) = ((q-1) - 2n + r)r$.  This codimension is greater than $2n$ when $r>1$ and $q \geq \frac72 n$, and so there is an open dense subset of $\hman$ consisting of $\tilde{g}$ with the property that $j^1\tilde{g}$ is disjoint from $S_r$ for $r > 1$, and that $(j^1\tilde{g})^{-1}(S_1)$ has codimension $q-2n$.

Now (see \cite{GolubitskyGuillemin}, Chapter VI, Lemma 3.2), a vector $v = (v_1,v_2)$ is tangent to $(j^1\tilde{g})^{-1}(S_1) = \Sigma_{\tilde{g}}$ if and only if $v$ is in the kernel of the intrinsic derivative
\begin{align}
\label{eqn:tangentvectorcondition}
D(d\tilde{g})_{(x,y)} : T_{(x,y)}(M \times M) \to \Hom(\Ker((d\tilde{g})_{(x,y)}),\Coker((d\tilde{g})_{(x,y)}))).
\end{align}
We would like to view this as a condition on the $2$-jet space.

Consider the projection $J^2(M \times M, \R^{q-1}) \to J^1(M \times M, \R^{q-1})$.  Let $S_1$ be as above, and let $S_1^2$ be its preimage in the $2$-jet space.  We have the following bundle map:
\[
\begin{tikzcd}
S_1^2 \arrow{drr}{} \arrow{rrr}{} &&& \Hom(T(M \times M), \Hom(\Ker((d\tilde{g})_{(x,y)}),\Coker((d\tilde{g})_{(x,y)}))) \arrow{dl}{} \\ 
&& S_1
\end{tikzcd}
\]
and the top arrow is a submersion.  In our case, $\Ker((d\tilde{g})_{(x,y)})$ is $1$-dimensional, so we may think of the upper-right space as $\Hom(T(M \times M), \Coker((d\tilde{g})_{(x,y)}))$.

We are now prepared to complete Step 2 in the outline of the proof of Proposition \ref{prop:main}.

\begin{lem} 
\label{lem:immersion}
For $q \geq \frac72 n$, the set $\mathscr{H}_{\mbox{\emph{imm}}} \cap \mathscr{H}_{\mbox{\emph{not}}}$ is open and dense in $\mathscr{H}_N^i$.
\end{lem}

\begin{proof}
We will define a certain subset $V$ of $J^2(M \times M, \R^{q-1})$ which describes the non-immersion condition, in the sense that $\tilde{g} \in \mathscr{H}_{\mbox{\emph{imm}}}$ if and only if its $2$-jet extension is disjoint from $V$.  The nowhere-tangent condition of $\hnot$ will be similar, because both conditions concern the tangent space of $\Sigma_{\tilde{g}}$.  Specifically, we can write
\begin{align*}
\tilde{g} \notin \himm & \Longleftrightarrow \mbox{ there exists } (x,y) \in \Sigma_{\tilde{g}} \mbox{ and } v_2 \in T_yM \mbox{ such that } (0,v_2) \in T_{(x,y)}\Sigma_{\tilde{g}}, \\
\tilde{g} \notin \hnot & \Longleftrightarrow \mbox{ there exists } (x,y) \in \Sigma_{\tilde{g}} \mbox{ and } v_2 \in T_yM \mbox{ such that } (s_x,v_2) \in T_{(x,y)}\Sigma_{\tilde{g}},
\end{align*}
where $s_x$ has been defined before Lemma \ref{lem:assumption}.

Incorporating this into (\ref{eqn:tangentvectorcondition}) and the subsequent diagram, a vector of the form $(0,v_2)$ is in $T_{(x,y)}\Sigma_{\tilde{g}}$ if and only if the map 
\[
T_yM \to T_xM \oplus T_yM \simeq T_{(x,y)}(M \times M) \to \Coker((d\tilde{g})_{(x,y)}).
\]
drops rank.  This is a linear map from an $n$-dimensional space to the $(q-1) - (2n - 1) = (q-2n)$-dimensional cokernel, and so the rank drop occurs on some submanifold $\Sigma$ of $T_{(x,y)}(M \times M) \to \Coker((d\tilde{g})_{(x,y)})$ with codimension $q - 3n + 1$.  As the top arrow in the diagram above is a submersion, the preimage $V$ of $\Sigma$ is a submanifold of $S_1^2$ with the same fiberwise codimension.  

Assume momentarily that the conditions of Transversality Theorem A apply to $V$, so that the set of $\tilde{g}$ such that $j^2(\tilde{g} \big|_{N^c}) \pitchfork V$ is open and dense.  For such $\tilde{g}$, $(j^2(\tilde{g} \big|_{N^c}))^{-1}(V)$ has codimension which is the sum of $q - 2n$ (the codimension in the $1$-jet fiber) and $q - 3n + 1$ (the codimension in the $2$-jet fiber).  The sum $2q-5n+1$ is greater than $2n$ when $q \geq \frac72 n$, so in these relative dimensions, transversality is disjointness.  Since $\himm$ is precisely the set of $\tilde{g}$ such that $j^2(\tilde{g} \big|_{N^c}) \cap V = \emptyset$, $\himm$ is open and dense.

The proof for $\hnot$ is similar, with the vector $(s_x,v_2)$ in place of $(0,v_2)$.  The only exception is that the map from $T_yM$ is not linear, but we can compute the codimension using the idea at the end of the proof of Theorem \ref{thm:strat}.  The pair $((s_x,v_2),(s_x,s_y))$ lies in the cokernel of $d\tilde{g}$ with codimension $q-2n$; taking the union over $v_2$ in the $(n-1)$-dimensional projectivization of $T_yM$ gives codimension $q - 3n +1$.

It remains to argue that $V$ is transverse to $\mathscr{J}$.  There is nothing to show on the $1$-jet level, as the $1$-jet portion of $\mathscr{J}$ is the entire $1$-jet.  Moreover, note that $V$ is independent of the source and target points, and that $V$ sits over only those points in the $1$-jet at which the rank drops by $1$.  Thus it is enough to consider fixed source and target points, a fixed linear map $P : \R^n \times \R^n \to \R^{q-1}$ of corank $1$, and a fixed symmetric bilinear map $B : (\R^n \times \R^n) \circ (\R^n \times \R^n) \to \R^{q-1}$, which we assume is an element of $V \cap \mathscr{J}$.  We would like to show transversality at $B$.

Let $s$ be a fixed nonzero vector in the (one-dimensional) kernel of $P$.  Since $B \in V$, there exists a vector $v = (0,v_2)$ such that $B(v, s)$ lies in the image of $P$.  Using rectangular coordinates orthogonal to the image of $P$, we may write $(B_1,\dots,B_{q-2n})$ for the projection to the cokernel of $P$.  Now since $B \in \mathscr{J}$, each component $B_i$ may be written, using the standard basis on $(\R^n \times \R^n) \circ (\R^n \times \R^n)$, in the following form:
\[
B_i = \left( \begin{array}{cc} A_i & 0 \\ 0 & D_i \end{array} \right).
\]
That is, the linear map $B_i(\cdot, s)$ can be represented as $(s_xA_i \ \ s_yD_i)(\cdot)$.  In this notation, the condition $B \in V$ is the condition that $s_yD_iv_2 = 0$ for all $i$.  We would like to show that the tangent space of $V$ at $B$ is transverse to $\mathscr{J}$.  As $B$ is an element of a space of linear maps, the tangent space at $B$ looks like the space itself, and the vectors tangent to $\mathscr{J}$, in the representation above, are those which have zeroes in the same spots as $B$ itself.  We will show transversality at $B$ by demonstrating a number of explicit linearly independent elements of the tangent space to $V$ which are transverse to $\mathscr{J}$.

For fixed $k,l$, $1 \leq k \leq 2n$ and $2n + 1 \leq l \leq 4n$, let $E_{kl}$ be the $4n \times 4n$ matrix with $\varepsilon$ in the $(k,l)$ and the $(l,k)$ spots.  

We compute $(B_i + E_{kl})(v,s) = s_yD_iv_2 + s_x^l v_2^k \varepsilon$; here the superscripts represent the $l$ and $k$ components of the respective vectors.  To arrange that this equals zero, we may perturb $A_i$ and $D_i$ appropriately to obtain a nearby map $B_i'$.  

Thus this new map, $(B_1,\dots,B_{i-1},B'_i + E_{kl},B_{i+1},\dots,B_{q-2n})$, obtained from an $\varepsilon$ movement transverse to $\mathscr{J}$ and a small movement (of $A_i$ or $D_i$) parallel to $\mathscr{J}$, still lies in $V$.  As this holds for any $i,k,l$, we have established transversality of $V$ to $\mathscr{J}$.
\end{proof}

\subsection{Step 3 of the proof of Proposition \ref{prop:main}: a local Taylor argument}
\label{sec:taylor}

We proceed with the proof of step 3, which resembles a certain classical argument involving immersions and embeddings.  For example, consider an immersion $a : X \to Y$, where $X$ is compact.  It is well known that for any point $x \in X$, there is a neighborhood of $x$ on which $a$ is an embedding. However, the following stronger statement holds:
there exists a neighborhood $U_x$ of $x$ and a neighborhood $W_a$ of $a$ in $C^\infty(X,Y)$ such that for any $b \in W_a$, $b|_{U_x}$ is an embedding.  Then by compactness of $X$, there exists $\delta$ (for example, the Lebesgue number for the cover $U_x$) such that any $b \in W_a$ is an embedding when restricted to any $\delta$-ball in $X$.

We require a similar statement for our purpose.  In particular, by definition of $\tilde{g} \in \himm$, the map $\pi |_{\Sigma_{\tilde{g}} \cap N^c}$ is an immersion.  Therefore, there exists a small number $\beta'$ such that $\pi|_{\Sigma_{\tilde{g}} \cap N^c}$ is injective on neighborhoods of size $\beta'$.  We claim that there exists an open neighborhood $W_{\tilde{g}}$ of $\tilde{g}$ in $\himm$ and a small number $\beta$, such that for any $\psi \in W_{\tilde{g}}$, the map $\pi|_{\Sigma_\psi \cap N^c}$ is injective on neighborhoods of size $\beta$.

The statement is essentially identical to the classical statement, except that the space $\Sigma_\psi$ is not generally the same as $\Sigma_{\tilde{g}}$; however, these submanifolds are close, so the projections are as well.

Using the same Lebesgue number argument as above, it is enough to establish a local version of the statement.

As usual, let $\Sigma_{\tilde{g}}$ be the set of points on which $d\tilde{g}$ has corank $1$, and let $S_1$ represent the submanifold of the $1$-jet consisting of linear maps with corank $1$, so that    $(x,y) \in \Sigma_{\tilde{g}}$ if and only if $j^1\tilde{g}(x,y) \in S_1$.  Fix such an $(x,y)$ and choose local coordinates $(z',z)$ at $j^1\tilde{g}(x,y)$ so that $z'$ parametrizes the submanifold $S_1$ near $(x,y)$ and $z$ parametrizes the normal space to $S_1$ at $(x,y)$. 

Shrink the domain of $\tilde{g}$ to a small neighborhood $D^n \times D^n$, call this restricted function $\varphi$.  We will assume that the neighborhood is small enough such that the image of $\varphi$, as well as the image of $\psi$ for any sufficiently small perturbation $\psi$ of $\varphi$, is contained within the scope of the local parametrization $(z',z)$.  

\begin{lem}
\label{lem:taylor}
Consider $\varphi \in C^\infty(D^n \times D^n, \R^{q-1})$, fix $(x_0,y_0) \in \Sigma_\varphi$, and let $(z',z)$ be the coordinates described above, localized at $j^1\varphi(x_0,y_0)$.  Assume that for the projection $\pi$ to the first factor, the restriction $\pi |_{\Sigma_\varphi}$ is an immersion.  Then there exists a neighborhood $W_{\varphi}$ of $\varphi$ in $C^\infty(D^n \times D^n, \R^{q-1})$ and a number $\eta > 0$ such that for all $\psi \in W_{\varphi}$, $\pi|_{\Sigma_\psi}$ is an embedding on every $\eta$-ball.
\end{lem}

\begin{proof}

For a point $(x,y) \in D^n \times D^n$ and a nearby point $(x,y+h)$, where $h = (h_1,\dots,h_n) \in \R^n$, we may use the Taylor formula to write
\begin{align*}
z_\psi(x,y+h) = z_\psi(x,y) + \frac{\partial z_\psi}{\partial y}(x,y)h + O(|h|^2);
\end{align*}
here $z_\psi$ represents the $z$-value of a function $\psi$.  By definition of the $(z',z)$ coordinates, $\Sigma_\psi$ is the set where $z_\psi = 0$.  Therefore, if we assume that $(x,y)$ and $(x,y+h)$ are both in $\Sigma_\psi$, it follows that
\begin{align}
\label{eqn:taylor}
\bigg\| \frac{\partial z_\psi}{\partial y}(x,y)h \bigg\| < C_\psi|h|^2
\end{align}
for some constant $C_\psi$.

The idea of the proof is to choose a sufficiently small neighborhood $W_\varphi$ of $\varphi$ such that this statement cannot hold for any $\psi \in W_\varphi$ and any points $(x,y)$, $(x,y+h)$ from $\Sigma_\psi$.  For this we need some estimates on $z_\varphi$ and its derivatives.

By the immersion hypothesis of the lemma, if $z_\varphi(x,y) = 0$, then the rank of $\frac{\partial z_\varphi}{\partial y}(x,y)$ is maximal.  It follows that there exists $\varepsilon > 0$ such that if $|z_{\varphi}(x,y)| \leq \varepsilon$, then the rank of $\frac{\partial z_\varphi}{\partial y}(x,y)$ is maximal.  For each such $(x,y)$, the norm of the operator $\frac{\partial z_\varphi}{\partial y}(x,y)$ is bounded below.  Then by compactness, there exists a global bound $\delta > 0$ such that 
\begin{align}
\label{eqn:onenorm}
\frac{ \| \frac{\partial z_\varphi}{\partial y}(x,y) h \|}{\| h \|} > 2\delta , \mbox{ for all } (x,y) \mbox{ such that } |z_{\varphi}(x,y)| \leq \varepsilon, \mbox{ and all nonzero } h \in \R^n.
\end{align}
Also by compactness of $D^n \times D^n$, the $2$-norm of $z_\varphi$ is bounded above by some constant; we write this constant in the form $C - \delta$ for some $C > \delta > 0$ (for the $\delta$ defined above).  

Let $W_\varphi$ be a neighborhood of $\varphi$ of size $\min\{\varepsilon,\delta\}$ in the $C^2$ norm.  Then for each $\psi \in W_\varphi$, we claim that the following hold:
\begin{enumerate}
\item For $(x,y) \in \Sigma_\psi$, $\| \frac{\partial z_\psi}{\partial y}(x,y) h \| > \delta \|h\|$.
\item The $2$-norm of $z_\psi$ is bounded above by $C$; in particular, we may replace $C_\psi$ in (\ref{eqn:taylor}) with the constant $C$, independent of $\psi$.
\end{enumerate}
The second item follows from the definition of $C$ and the size of $W_\varphi$. For the first, we have $(x,y) \in \Sigma_\psi$, hence $z_\psi(x,y) = 0$, hence $|z_\varphi(x,y)| \leq \varepsilon$ by the choice of $W_\varphi$.  Therefore (\ref{eqn:onenorm}) holds, hence the first item holds by the choice of $W_\varphi$.

Finally, let $\eta = \frac{\delta}{C}$.  If $(x,y)$ and $(x,y+h)$ are both contained in an $\eta$-ball in $\Sigma_\psi$, then $|h| < \eta$, which in conjunction with item 1, yields that the quantity in (\ref{eqn:taylor}) is $> C\|h\|^2$, contradicting item 2.
\end{proof}

\subsection{Applications of Transversality Theorem B to Step 4}
\label{sec:ttb}

We now carry out Step 4 of the proof of Proposition \ref{prop:main}.  Let us first consider the condition that there are no far-away double points of $\pi|_{\Sigma_{\tilde{g}}}$; this is a slightly stronger condition than required for Step 4.  Let $\alpha$, $\beta$ be as defined in Steps 1 and 3 of the proof of Proposition \ref{prop:main}, and let $\gamma = \min\left\{\alpha, \beta\right\}$.  Let $N_\gamma$ be the $\gamma$-neighborhood of the diagonal of $\Delta M$, and we study the condition that there do not exist $(x,y)$, $(x,y')$ both in $\Sigma_{\tilde{g}} \cap N_\gamma^c$ which are of distance at least $2\beta$ apart.  We may formulate this as a condition on $X = M \times M - N_\gamma$.  Let $Z_\beta$ be the $\beta$-neighborhood of the diagonal $\Delta X$.  Then $X \times X - Z_\beta$ consists of pairs of points of $M$ which are at least $\gamma$ away from the diagonal and $2\beta$ apart.  

By definition of $\Sigma_{\tilde{g}}$ we are concerned with points $((x,y),(x,y')) \in X \times X - Z_\beta$ with a simultaneous rank drop of $d{\tilde{g}}$ at $(x,y)$ and $(x,y')$.  Informally, each rank drop occurs with codimension $q - 2n$, so both happen simultaneously with codimension $2q - 4n$, which is larger than $3n$ provided that $q > \frac72 n$.  (Note that this is the same informal bound one achieves by dimension-counting the condition that $\pi\big|_{\Sigma_{\tilde{g}}}$ is an embedding.) 

More formally, let $S_r \subset J^1(X, \R^{q-1})$ be the set which consists fiberwise of linear maps with corank $r$, and let $\Sigma \subset X \times X - Z_\beta$ be the subset consisting of points of the form $(x,y,x,y')$.  Let $V$ be the subset of $J_2^1(X, \R^{q-1})$ which fiberwise is $S_r \times S_r$ over points of $\Sigma$ and is empty over points of $\Sigma^c$.  Then $V$ has codimension $2q-4n$, and Transversality Theorem B applies (since in the $2$-fold-$1$-multijet, $\mathscr{J}_2$ is the entire space and there is no transversality of $V$ to show).  Therefore, the space of $\tilde{g} \in \mathscr{H}_{N_\gamma}^i(\tilde{h})$ such that $j_2^1(\tilde{g})|_{Z_\beta^c} \pitchfork V$ is open and dense.  For such $\tilde{g}$, $(j_2^1(\tilde{g} \big|_{\Sigma}))^{-1}(V)$ has codimension $2q-4n$, hence is empty when $q > \frac72 n$, and assumptions 4 and 5 of Lemma \ref{lem:assumption} are guaranteed for generic $\tilde{g}$.

In case $q = \frac72 n$, $(j_2^1(\tilde{g} \big|_{\Sigma}))^{-1}(V)$ is generically a finite set; consisting of points of the form $(x,y), (x,y')$ which are both in $\Sigma_{\tilde{g}}$.  We would like to ensure that the projected sections $s_x(x,y)$ and $s_x(x,y')$ are not parallel.  This can be achieved by a small modification of the above argument.  We may define $V'$ as the subset of $S_1 \times S_1$ consisting of pairs of linear maps $P,P' \in \Hom(\R^n \times \R^n, \R^{q-1})$ both of corank $1$, such that $\ker P = \ker P'$.  This set has codimension $q - 2n$ and is missed by the $2$-fold $1$-multijet of $\tilde{g}$ for generic $\tilde{g}$ provided that $q > \frac52 n$.

To study the no-triple-point condition, one may state and prove the obvious $3$-multijet version of Transversality Theorem B, from which it follows that triple points of $\pi \big|_{\Sigma_{\tilde{g}}}$ do not exist for generic $\tilde{g}$ as long as $q > \frac{10}{3} n$.

\section{Statements of transversality}
\label{sec:transversality}

We begin with a transversality theorem due to Haefliger, similar to the multijet transversality theorem.

\begin{thm}[Haefliger \cite{Haefliger}, 1.8--1.10] Let $M$ be a compact manifold, $V \subset J^s(M, \R^{i-1}) \times J^s(M, \R^{i-1}) $ be a stratified set, and let $N$ be a fixed open neighborhood of $\Delta M$.  Then 
\begin{align*}
\left\{ F \in C^\infty(M, \R^{i-1}) \ \Big| \ j^s_2F|_{N^c} \pitchfork V \right\}
\end{align*}
is open and dense in $C^\infty(M,\R^{i-1})$.
\label{thm:haefliger}
\end{thm}

We would like to formulate a similar statement which will apply to the functions we study here.  We would like our statement to apply to functions of the form $F \oplus F : M \times M \to \R^{i-1}$.  In this case, $j^s(F \oplus F)$ is a section of the bundle $J^s(M \times M, \R^{i-1})$, yet we cannot expect that for an arbitrary $V$, a generic $F$ will satisfy $j^s(F \oplus F) \pitchfork V$.  Indeed, such $s$-jets do not occupy the entire jet space $J^s(M \times M, \R^{i-1})$, and so it is easy to choose a $V$ for which this statement fails.  (To see this, compare the possible forms of sections $j^2(F \oplus F)$ with sections $j^2\varphi$ for arbitrary functions $\varphi : M \times M \to \R^{i-1}$).

Thus, we let $\mathscr{J}'$ be the image of the natural inclusion $J^s(M, \R^{i-1}) \times J^s(M, \R^{i-1}) \to J^s(M \times M, \R^{i-1})$, and we claim that our proposed statement is true for those $V$ satisfying $V \pitchfork \mathscr{J}'$.

\begin{cor}
\label{cor:haefliger}
Let $M$ be a compact manifold, $V \subset J^s(M \times M, \R^{i-1})$ be a stratified set satisfying $V \pitchfork \mathscr{J}'$, and let $N$ be a fixed open neighborhood of $\Delta M$.  Then 
\begin{align*}
\left\{ F \in C^\infty(M, \R^{i-1}) \ \Big| \ j^s(F \oplus F)|_{N^c} \pitchfork V \right\}
\end{align*}
is open and dense in $C^\infty(M,\R^{i-1})$.
\end{cor}

\begin{proof} The functions $F$ satisfying the upshot of the corollary for some $V \subset J^s(M \times M, \R^{i-1})$ are precisely those satisfying Haefliger's theorem applied to $V \cap \mathscr{J}'$. 
\end{proof}

For the remainder of this section, $N$ will represent a fixed open neighborhood of the diagonal, and let $N' \subset N$ be a smaller open $\Z_2$-invariant neighborhood of the diagonal.  Consider the function space
\[
C \coloneqq C^\infty(M,\R^{i-1}) \times C^\infty((M \times M - N')/\Z_2,\R^{q-i})
\]
and let $\mathscr{J}$ be the subbundle of $J^s((M \times M - N')/\Z_2, \R^{q-1})$ generated by $s$-jets of functions of the form $(F \oplus F \big|_{(M \times M - N')/\Z_2}, G)$ for $(F,G) \in C$.  That is, $\mathscr{J}$ is the image of the inclusion
\[
J^s(M, \R^{i-1}) \times J^s(M, \R^{i-1}) \times J^s((M \times M - N')/\Z_2, \R^{q-i}) \to J^s((M \times M - N')/\Z_2, \R^{q-1}).
\]

Combining Corollary \ref{cor:haefliger} for $F$ with ordinary Thom transversality for $G$ yields the following.

\begin{cor}
\label{cor:mainlemmaA}
Let $M$ be a compact manifold, $N' \subset N$ fixed open neighborhoods of the diagonal $\Delta M$, and $V \subset J^s((M \times M - N')/\Z_2, \R^{q-1})$ be a stratified set such that $V \pitchfork \mathscr{J}$.  Then 
\begin{align*}
\left\{ (F,G) \in C \ \Big| \ j^s((F \oplus F) \big|_{(M \times M - N')/\Z_2},G) |_{N^c} \pitchfork V \right\}
\end{align*}
is open and dense in $C$.
\end{cor}

\begin{thmn}[Transversality Theorem A] Let $M$ be a compact manifold, $N' \subset N$ fixed open neighborhoods of the diagonal $\Delta M$, and $V \subset J^s((M \times M - N')/\Z_2, \R^{q-1})$ be a stratified set such that $V \pitchfork \mathscr{J}$.  Then
\begin{align*}
\left\{ \tilde{g} = (F \oplus F, G) \in \mathscr{H}_N^i(\tilde{h}) \ \Big| \ j^s(\tilde{g}) |_{(M \times M) \setminus N} \pitchfork V \right\}
\end{align*}
is open and dense in $\mathscr{H}_N^i(\tilde{h})$.
\end{thmn}

\begin{remark}
We briefly address a subtle distinction between the statement above and the original statement and applications of Transversality Theorem A.  For the purposes of verifying the transversality theorems, it is easier to consider $\tilde{g} : (M \times M - N')/ \Z_2 \to \R^{q-1}$, since then $\Z_2$-equivariance is built into the domain, and the domain is still a manifold.  In any case, since the transversality theorems apply only to $\tilde{g} \big|_{N^c}$, the actual domain is somewhat moot.
\end{remark}

\begin{proof}
We will consider $N$, $i$, and $\tilde{h}$ fixed, and write $\mathscr{H}$ for $\mathscr{H}_N^i(\tilde{h})$.  For $\delta \geq 0$, let $N_\delta$ refer to the open $\delta$-neighborhood of $N$, and apply the previous lemma to each $N_\delta$ to obtain that
$$
\Omega_\delta \coloneqq \left\{ (F,G) \in C \ \Big| \ j^s(F \oplus F, G) |_{N_\delta^c} \pitchfork V \right\}
$$
is open and dense in $C$.

Our goal is to show that $\Omega_0 \cap \mathscr{H}$ is open and dense in $\mathscr{H}$.  Consider the inclusion map $\iota : \mathscr{H} \to C$, then $\Omega_\delta \cap \mathscr{H} = \iota^{-1}(\Omega_\delta)$ is open.  In particular, $\Omega_0 \cap \mathscr{H}$ is open, and so it only remains to show that $\Omega_0 \cap \mathscr{H}$ is dense in $\mathscr{H}$.

Let $(F,G) \in \mathscr{H} \subset C$.  By density of $\Omega_\delta$ in $C$, there exists an approximation $(F_n,G_n)$ of $(F,G)$, where each $(F_n,G_n)$ is in $\Omega_\delta$ for some fixed $\delta$.  Now define a new function which is equal to $(F_n \oplus F_n,G_n)$ on $N_\delta^c$, is equal to $\tilde{h} \oplus \tilde{h}$ on $N$ (and hence its closure), and smoothly extends to $N_\delta \setminus N$, sufficiently approximating $(F,G)$ there.  This new sequence converges to $(F,G)$ and is in $\Omega_\delta \cap \mathscr{H}$.

Thus the set $\Omega_\delta \cap \mathscr{H}$ is open and dense in $\mathscr{H}$, so by the Baire category theorem, 
\[
\Omega_0 \cap \mathscr{H} = \bigcap_{m \in \N} \Omega_\frac{1}{m} \cap \mathscr{H}
\]
is dense in $\mathscr{H}$.
\end{proof}

Transversality Theorem B may be proven in much the same manner as Transversality Theorem A.  Let $M$ be a compact manifold and $N$ some open neighborhood $\Z_2$-invariant neighborhood of the diagonal.  Let $X = M \times M - N$, and let $\mathscr{J}_2^\prime = \mathscr{J}'|_X \times \mathscr{J}'_X$, where each factor represents the restriction of $\mathscr{J}' \subset J^s(M \times M, \R^{i-1})$ to the base space $X$.  The following corollary may be derived from Haefliger's Theorem applied to $X$, similar to our derivation of Corollary \ref{cor:haefliger} from Haefliger's Theorem above.

\begin{cor}
\label{cor:ttb}
Let $M$ be a compact manifold, $N'$ a fixed open $\Z_2$-invariant neighborhood of the diagonal $\Delta M$, $X = M \times M - N'$, $Z$ a fixed open neighborhood of the diagonal $\Delta X$, and $V \subset J_2^s(X,\R^{i-1})$ a stratified set such that $V \pitchfork\mathscr{J}_2^\prime$.  Then
\[
\left\{ F \in C^{\infty}(M, \R^{i-1}) \ \big| \ j_2^s\big( (F \oplus F)|_X\big)\big|_{Z^c} \pitchfork V \right\}
\]
is open and dense.
\end{cor}

Now let $C$ be as above and $\mathscr{J}_2 = \mathscr{J}|_{X/\Z_2} \times \mathscr{J}_{X / \Z_2}$.

\begin{cor}
\label{cor:ttb2}
Let $M$ be a compact manifold, $N'$ a fixed open $\Z_2$-invariant neighborhood of the diagonal $\Delta M$, $X = M \times M - N'$, $Z$ a fixed open neighborhood of the diagonal $\Delta X$, and $V \subset J_2^s(X/\Z_2,\R^{q-1})$ a stratified set such that $V \pitchfork\mathscr{J}_2$.  Then
\[
\left\{ (F,G) \in C \ \big| \ j_2^s\big( (F \oplus F)|_{X/\Z_2},G\big)\big|_{Z^c} \pitchfork V \right\}
\]
is open and dense.
\end{cor}

Transversality Theorem B may then be proven in the same manner as Transversality Theorem A.

We state one final variant, obtained by assuming no neighborhood $Z$, and considering the special case where $i = q = 4n$ and $s=0$.  With the latter assumption, $\mathscr{J}_2$ is the entire $2$-fold $0$-multijet and so there is no ``transversality to $\mathscr{J}_2$'' assumption.  This variant is used exclusively in the proof of Proposition \ref{prop:admissibleexist}.

\begin{cor}
\label{cor:ttb3}
Let $M$ be a compact manifold, $N'$ a fixed open $\Z_2$-invariant neighborhood of the diagonal $\Delta M$, $X = M \times M - N'$, and $V \subset J_2^0(X/\Z_2,\R^{4n-1})$ a stratified set.  Then
\[
\left\{ h \in C^\infty(M,\R^{4n-1}) \ \big| \ j_2^0\big( (h \oplus h)|_{X/\Z_2}) \pitchfork V \right\}
\]
is a residual set.
\end{cor}

\section{Generalizations and Conjectures}
\label{sec:conclusion}

The assumptions of Lemma \ref{lem:assumption} are shown to hold using genericity arguments, but it may be possible to instead use some kind of h-principle statements for $\pi$.  For example, to show that $\pi$ is generically an immersion, we compute the expected dimension of $\Sigma_{\tilde{g}}$ is $4n-q$, and then use the fact that a $(4n-q)$-dimensional space generically immerses in an $n$-dimensional space as long as $2(4n-q) \leq n$, or $q \geq \frac72 n$.  However, it is possible that $\pi$ is an immersion even when the dimension of $\Sigma_{\tilde{g}}$ is exactly $n$, i.e.\ $q = 3n$.  Perhaps it is possible to show that, instead of making perturbations of $h$ and $g$ to achieve the assumptions of Lemma \ref{lem:assumption}, we may homotope them to maps which satisfy the assumptions.

In any case, we believe that the h-principle-type result of Conjecture \ref{con:main} (in Section \ref{sec:intro}) holds for totally nonparallel immersions.  This conjecture has the flavor of Haefliger's theorem on embeddings \cite{Haefliger}, which states that for $q > \frac32 n + \frac32$, the existence of an embedding of a compact $n$-dimensional manifold $M$ to a $q$-dimensional manifold $N$ is equivalent to the existence of a map $g : M \times M \to_{\Z_2} N \times N$ such that $g^{-1}(\Delta N) = \Delta M$.  In this case, the $\Z_2$-equivariant map near the diagonal may be realized as a \emph{skew bundle map} \cite{HaefligerHirsch} (unrelated to skew loops and totally skew embeddings) studied by Haefliger and Hirsch, who showed that in the appropriate relative dimensions, the space of skew bundle maps is homotopy equivalent to the space of bundle monomorphisms.  Then, one may apply the Smale-Hirsch theorem to homotope a monomorphism to an immersion.  Thus the local problem is solved, and the removal of singularities technique can be applied, outside the diagonal, in much the same manner as we have done here.

In our case, we would need to show that the immersion of the unordered configuration space, near the diagonal, is homotopic (in the appropriate function space) to a semifree map.  This would likely require an algebraic Haefliger-Hirsch-type statement and then an h-principle for semifree maps, each of which has eluded us even after a considerable amount of effort.

We have also made some attempts to study the totally skew embedding condition defined in the Introduction. The difficulty in importing our results to the totally skew situation arises due to the complexity of the local condition.  A sufficient local condition for total skewness for curves $\R \to \R^3$ is nonzero curvature and nonzero torsion.  We were unable to generalize this to find the correct local condition for higher-dimensional total skewness, and so we were not able to repeat Section \ref{sec:local} for totally skew maps.  However, the results of the later sections apply, a version of Lemma \ref{lem:assumption} can be formulated, and we can still perform the inductive replacement provided the local problem is solved.

We conclude this paper with one final example to distinguish between semifreedom and the totally nonparallel condition.  Perhaps a generalization of this example could yield a totally nonparallel immersion $\R^3 \to \R^7$.

\subsection{A low-dimensional example and the hyperdeterminant}
\label{sec:low}

We conclude with a short discussion of totally nonparallel and semifree maps $\R^2 \to \R^4$.  In particular, we give some examples to distinguish the totally nonparallel condition from the semifree condition, to give hope that totally nonparallel immersions might exist even in dimensions whether the corresponding nonsingular symmetric bilinear maps do not exist.

To begin, observe that for a smooth function $f : \R \to \R$, the condition $f'' \neq 0$ is sufficient for the graph of $f$ to be totally nonparallel.  It is not necessary; take $f(x) = x^4$.  Of course, if the condition $f'' = 0$ holds on a neighborhood, the graph of $f$ cannot be totally nonparallel.

Now consider a quadratic map $\R^2 \to \R^2 : (x,y) \mapsto (f(x,y),g(x,y))$ whose graph is totally nonparallel.  By Theorem \ref{thm:snbmap}, the associated symmetric bilinear map $B$ is nonsingular.  Here $B$ represents the (constant) Hessian of $(f,g)$.  Thus we have, for all $a \neq b$ and all $u \neq 0$,
$$\left[ (a_1 - b_1) \left( \begin{array}{cc} f_{xx} & f_{yx} \\
g_{xx}& g_{yx} \end{array} \right) + (a_2-b_2) \left( \begin{array}{cc} f_{xy} & f_{yy} \\
g_{xy}& g_{yy} \end{array} \right) \right] \left( \begin{array}{c} u_1 \\ u_2 \end{array} \right) \neq 0.$$
Therefore the linear combination is an invertible matrix for all $a \neq b$, and so the following determinant of the linear combination is nonzero for all $a \neq b$:
\begin{align*}
(f_{xx}g_{xy}-f_{xy}g_{xx})(a_1-b_1)^2 & + (f_{xx}g_{yy}-f_{yy}g_{xx})(a_1-b_1)(a_2-b_2) \\ & + (f_{xy}g_{yy}-f_{yy}g_{xy})(a_2-b_2)^2.
\end{align*}
Hence the following discriminant is negative:
\begin{align}
\label{eqn:hyper}
0 > (f_{xx}g_{yy}-f_{yy}g_{xx})^2 - 4(f_{xx}g_{xy}-f_{xy}g_{xx})(f_{xy}g_{yy}-f_{yy}g_{xy}).
\end{align}
The quantity in the right side of (\ref{eqn:hyper}) is Cayley's \emph{hyperdeterminant} (see \cite{GelfandKapranovZelevinsky}) of the Hessian of $(f,g)$, considered as a $2 \times 2 \times 2$ array.

It is tempting to guess that \emph{any} totally nonparallel immersion must have negative hyperdeterminant at some point, and therefore induces a symmetric nonsingular bilinear map.  A result of this nature might indicate that the problem of totally nonparallel immersions is equivalent to the problem of symmetric nonsingular bilinear maps.  But the following example suggests otherwise.

\begin{example}
The graph of the function $(x,y) \mapsto (x^4-y^4, x^3y+xy^3)$ is a totally nonparallel immersion of $\R^2$ into $\R^4$ with hyperdeterminant identically $0$.  To check the first claim, note that the graph of a function $\R^n \to \R^n$ is totally nonparallel if and only if its Jacobian $J$ satisfies $\det(J(a) - J(b)) \neq 0$ for all $a \neq b \in \R^n$.  Here we compute
\begin{align*}
0 < \operatorname{det}(J(a)-J(b)) = 
& (a_1^2+a_1b_1 + b_1^2 + a_2^2+a_2b_2+b_2^2) \cdot  \\ & \left[ ((a_1-b_1)^2+(a_2-b_2)^2)^2+3(a_1^2+a_2^2-b_1^2-b_2^2)^2 \right] ,
\end{align*}
so that the graph is totally nonparallel.

On the other hand, if we write $f(x,y) = x^4-y^4$ and $g(x,y) = x^3y+xy^3$ we have
\begin{align*}
f_{xx} & = 12x^2, \hspace{.94in} f_{xy} =0, \hspace{1.25in} f_{yy} = -12y^2, \\
g_{xx} & = 6xy, \hspace{1in} g_{xy} =3x^2+3y^2, \hspace{.70in} g_{yy} = 6xy,
\end{align*}
so the hyperdeterminant of the Hessian is
\[
(72x^3y+72xy^3)^2 - 4(36x^4+36x^2y^2)(36x^2y^2+36y^4) = 0.
\]
\end{example}

We know by Corollary \ref{cor:noquad} that there is no quadratic totally nonparallel immersion $\R^3 \to \R^7$.  Can the above example be modified to produce a non-quadratic totally nonparallel immersion $\R^3 \to \R^7$?

\bibliographystyle{plain}
\bibliography{../../../bib}{}

\end{document}